\documentclass[aap,MSNbibl,seceqn,nameyear,dvips]{arximspdf}

\doi{10.1214/11-AAP835} 
\volume{22}
\issue{6}
\pubyear{2012}
\firstpage{2388}
\lastpage{2428}

\makeatletter

\newtheorem{theorem}{Theorem}[section]
\newtheorem{proposition}[theorem]{Proposition}
\newtheorem{corollary}[theorem]{Corollary}
\newtheorem{lemma}[theorem]{Lemma}

\newproclaim{example}[theorem]{Example}

\newproclaim{Rems}{Remarks}
\newproclaim{Rem}{Remark}

\newcommand{\Cov}{\operatorname{Cov}}
\newcommand{\dd}{\,d}
\newcommand{\p}{{\bolds{\mathsf p}}}

\newcommand{\essinf}{\mathop{\operatorname{ess}\operatorname{inf}}}

\newcommand{\dpP}[1]{#1^{\p}}
\newcommand{\sdpP}[1]{#1^{\p}}
\newcommand{\spP}[1]{{}^{\p}#1}

\newcommand{\sSIdot}{\cdot}
\newcommand{\SIdot}{\cdot}

\newcommand{\Var}{\operatorname{Var}}

\newcommand{\Lt}{{\widetilde L}}
\newcommand{\Qt}{\widetilde Q}
\newcommand{\thetat}{\widetilde\vartheta}
\newcommand{\Vt}{\widetilde V}
\newcommand{\Qh}{\widehat Q}
\newcommand{\Wh}{\widehat W}
\newcommand{\Zh}{\widehat Z}
\newcommand{\A}{{\mathcal A}}
\newcommand{\B}{{\mathcal B}}
\newcommand{\E}{{\mathcal E}}
\newcommand{\F}{{\mathcal F}}
\newcommand{\G}{{\mathcal G}}
\renewcommand{\H}{{\mathcal H}}

\newcommand{\M}{{\mathcal M}}
\newcommand{\N}{{\mathcal N}}
\renewcommand{\P}{{\mathcal P}}
\renewcommand{\SS}{{\mathcal S}}

\newcommand{\FF}{{\mathbb F}}
\newcommand{\MM}{\langle M\rangle}
\newcommand{\MMc}{\langle M^c\rangle}
\newcommand{\MMd}{\langle M^d\rangle}

\newcommand{\NN}{{\mathbb N}}
\newcommand{\PP}{{\mathbb P}}
\newcommand{\RR}{{\mathbb R}}

\newcommand{\lsb}{[\![}
\newcommand{\rsb}{]\!]}

\newcommand{\anull}{a^{(0)}}
\newcommand{\aone}{a^{(1)}}
\newcommand{\Bq}{B^{\q}}
\newcommand{\deltac}{\delta^c}
\newcommand{\deltad}{\delta^d}
\newcommand{\fatdot}{{\cdot}}
\newcommand{\func}{g}
\newcommand{\Ls}{L^\prime}
\newcommand{\Lone}{L^{(1)}}
\newcommand{\Ltwo}{L^{(2)}}
\newcommand{\loc}{{\mathrm{loc}}}
\newcommand{\mnull}{m^{(0)}}
\newcommand{\mone}{m^{(1)}}
\newcommand{\Mc}{M^c}
\newcommand{\Md}{M^d}
\newcommand{\Nq}{N^{\q}}
\newcommand{\NY}{N^Y}
\newcommand{\NNq}[1]{{\mathcal N}_{#1}(\q)}
\newcommand{\NNY}[1]{{\mathcal N}_{#1}(Y)}
\newcommand{\Nnull}{N^{(0)}}

\newcommand{\pq}{{}^{\p}\q}
\newcommand{\pY}{{}^{\p}Y}
\newcommand{\PPes}{\PP_{\mathrm{e}, \sigma}^2(S)}
\newcommand{\PPss}{\PP_{\mathrm{s}, \sigma}^2(S)}
\newcommand{\psione}{\psi^{(1)}}
\newcommand{\psitwo}{\psi^{(2)}}
\newcommand{\phione}{\varphi^{(1)}}
\newcommand{\phitwo}{\varphi^{(2)}}
\newcommand{\psic}{\psi^c}
\newcommand{\psid}{\psi^d}

\newcommand{\q}{q}

\newcommand{\tE}{{}^t{\mathcal E}}
\newcommand{\thetast}{\vartheta^{*,t}}
\newcommand{\thetasn}{\vartheta^{*,0}}
\newcommand{\uE}{{}^u{\mathcal E}}

\newcommand{\VH}{V^H}
\newcommand{\Vn}{V^0}
\newcommand{\vi}{v^{(i)}}
\newcommand{\vnull}{v^{(0)}}
\newcommand{\vone}{v^{(1)}}
\newcommand{\vtwo}{v^{(2)}}

\newcommand{\Wb}{\bar W}
\newcommand{\Wone}{W^{(1)}}
\newcommand{\Wtwo}{W^{(2)}}

\newcommand{\Ynull}{Y^{(0)}}
\newcommand{\Yone}{Y^{(1)}}
\newcommand{\Ytwo}{Y^{(2)}}

\newcommand{\ZQ}{Z^Q}
\newcommand{\ZQP}{Z^{Q;P}}
\newcommand{\ZRP}{Z^{R;P}}
\newcommand{\ZRQ}{Z^{R;Q}}

\newcommand{\tr}{^{\mathrm{tr}}}

\makeatother

\begin{document}
\begin{frontmatter}

\title{Mean-variance hedging via stochastic control and BSDEs for
general semimartingales}
\runtitle{Mean-variance hedging and BSDEs for semimartingales}

\begin{aug}
\author[C]{\fnms{Monique} \snm{Jeanblanc}\thanksref{t1}\ead[label=e1]{monique.jeanblanc@univ-evry.fr}},
\author[A]{\fnms{Michael} \snm{Mania}\thanksref{t1,t2}\ead[label=e2]{misha.mania@gmail.com}},
\author[B]{\fnms{Marina}~\snm{Santacroce}\thanksref{t3}\ead[label=e3]{marina.santacroce@polito.it}}
\and
\author[D]{\fnms{Martin} \snm{Schweizer}\corref{}\thanksref{t4}\ead[label=e4]{martin.schweizer@math.ethz.ch}}
\runauthor{Jeanblanc, Mania, Santacroce and Schweizer}
\affiliation{Universit\'e d'Evry Val d'Essonne, Georgian-American
University,
Politecnico~di~Torino and ETH Z\"urich}
\address[A]{M. Mania\\
A. Razmadze Mathematical Institute\\
Aleksidze St. 1, Tbilisi 0193\\
and\\
Georgian-American University\\
17a, Chavchavadze Ave\\
Tbilisi\\
Georgia\\
\printead{e2}}
\address[B]{M. Santacroce\\
Dipartimento di Matematica\\
Politecnico di Torino\\
Corso Duca degli Abruzzi 24\\
I-10129 Torino\\
Italy\\
\printead{e3}}
\address[C]{M. Jeanblanc\\
Laboratoire Analyse et Probabilit\'es\\
Universit\'e d' Evry Val d' Essonne\\
IBGBI, 23 Boulevard de France\\
F-91037 Evry Cedex\\
France\\
and\\
Institut Europlace de Finance\\
Palais Brongniart\\
Rue de la Bourse\\
F-75002 Paris\\
France\\
\printead{e1}} 
\address[D]{M. Schweizer\\
Departement Mathematik\\
ETH Z\"urich\\
ETH-Zentrum, HG G 51.2\\
R\"amistrasse 101\\
CH-8092 Z\"urich\\
Switzerland\\
and\\
Swiss Finance Institute\\
Walchestrasse 9\\
CH-8006 Z\"urich\\
Switzerland\\
\printead{e4}\hspace*{14.1pt}}
\end{aug}

\thankstext{t1}{Supported by Chaire de Cr\'edit, French Banking
Federation.}

\thankstext{t2}{Supported by the Georgian National
Foundation Grant N. ST09\_471\_3-104.}

\thankstext{t3}{Supported in part by MIUR, Grant 2008YYYBE4.}

\thankstext{t4}{Supported by the Swiss
Finance Institute (SFI) and the National Centre of Competence in
Research ``Financial Valuation and Risk Management'' (NCCR
FINRISK), Project D1 (Mathematical Methods in Financial Risk
Management). The NCCR FINRISK is a research instrument of the Swiss
National Science Foundation.}

\received{\smonth{2} \syear{2011}}
\revised{\smonth{11} \syear{2011}}

%
\begin{abstract}
We solve the problem of mean-variance hedging for general
semimartingale models via stochastic control methods. After proving
that the value process of the associated stochastic control problem has
a quadratic structure, we characterize its three coefficient processes
as solutions of semimartingale backward stochastic differential
equations and show how they can be used to describe the optimal trading
strategy for each conditional mean-variance hedging problem. For
comparison with the existing literature, we provide alternative
equivalent versions of the BSDEs and present a number of simple
examples.
\end{abstract}

%
\begin{keyword}[class=AMS]
\kwd{60G48}
\kwd{60H10}
\kwd{93E20}
\kwd{91G80}.
\end{keyword}
\begin{keyword}
\kwd{Mean-variance hedging}
\kwd{stochastic control}
\kwd{backward stochastic differential equations}
\kwd{semimartingales}
\kwd{mathematical finance}
\kwd{variance-optimal martingale measure}.
\end{keyword}

\end{frontmatter}

\setcounter{section}{-1}
\section{Introduction}\label{sec0}

Mean-variance hedging is one of the classical problems from
mathematical finance. In financial terms, its goal is to minimize the
mean squared error between a given payoff $H$ and the final wealth of a
self-financing strategy $\vartheta$ trading in the underlying assets $S$.
Mathematically, one wants to project the random variable $H$ in
$L^2(P)$ on the space
of all stochastic integrals $\vartheta\SIdot S_T =
\int_0^T \vartheta_r \dd S_r$, perhaps after subtracting an initial
capital~$x$. The contribution of our paper is to solve this problem via
stochastic control methods and stochastic calculus techniques for the
case where the asset prices $S$ are given by a general (locally
$P$-square-integrable) semimartingale, under a natural no-arbitrage
assumption.

The literature on mean-variance hedging is vast, and we do not try to
survey it here; see \citet{Sch} for an attempt in that direction.
There are
two main approaches; one of them uses martingale theory and projection
arguments, while the other views the task as a linear-quadratic
stochastic control problem and uses backward stochastic differential
equations (BSDEs) to describe the solution. By combining tools from
both areas, we improve earlier work in two directions---we describe
the solution more explicitly than by the martingale and projection
method, and we work in a general semimartingale model without
restricting ourselves to particular setups (like It\^o processes or
L\'evy settings). We show that the value process of the stochastic control
problem associated to mean-variance hedging possesses a quadratic
structure, describe its three coefficient processes by
semimartingale BSDEs and show how to obtain the optimal strategy
$\vartheta^*$ from there. In contrast to the majority of earlier
contributions from the control strand of the literature, we also give a
rigorous derivation of these BSDEs. For comparison, the usual results
(especially in settings with It\^o processes or jump-diffusions) start
from a BSDE system and only prove a verification theorem that shows how
a solution to the BSDE system induces an optimal strategy. Apart from
being more precise, we think that our approach is also more informative
since it shows clearly and explicitly how the BSDEs arise, and hence
provides a systematic way to tackle mean-variance hedging via
stochastic control in general semimartingale models. More detailed
comparisons to the literature are given in the respective sections.

The paper is structured as follows. We start in Section~\ref{sec1} with a
precise problem formulation and state the martingale optimality
principle for the value process $\VH(x)$ of the associated stochastic
control problem. Assuming that each (time $t$) conditional problem
admits an optimal strategy, we then show that $\VH(x)$ is a quadratic
polynomial in $x$ whose coefficients are stochastic processes $\vnull,
\vone, \vtwo$ that do not depend on $x$. This is a kind of folklore
result, and our only claim to originality is that we give a very simple
proof in a very general setting. We also show that the coefficient
$\vtwo$ equals the value process $\Vn(1)$ for the control problem with
initial value $x=1$ and $H\equiv0$.

Motivated by the last result, we study in Section~\ref{sec2} the
particular problem for \mbox{$x=1$} and $H\equiv0$. We impose the no-arbitrage
condition that there exists an equivalent $\sigma$-martingale measure
for $S$ with $P$-square-integrable density and are then able to
characterize the process $\vtwo$ as the solution of a semimartingale
BSDE. More precisely, Theorem~\ref{Theorem2.4} shows that all conditional
problems for $x=1, H\equiv0$ admit optimal strategies if and only if
that BSDE (\ref{equ2.18}) has a solution in a specific class, and in that
case, the unique solution is $\vtwo$ and the conditionally optimal
strategies can be given in terms of the solution to (\ref{equ2.18}).\vadjust{\goodbreak} In
comparison to earlier work, we eliminate all technical assumptions
(like continuity or quasi-left-continuity) on $S$, and we also do not
need reverse H\"older inequalities for our main results.

Section~\ref{sec3} considers the general case of the mean-variance
hedging problem with $x\in\RR$ and $H\in L^2(\F_T,P)$. The analog of
Theorem~\ref{Theorem2.4} is given in Theorem~\ref{Theorem3.1}, where we describe
the three coefficient processes $\vtwo, \vone, \vnull$ by a coupled
system (\ref{equ3.1})--(\ref{equ3.3}) of semimartingale BSDEs. Existence
of optimal strategies for all conditional problems for $(x,H)$ is shown
to be equivalent to solvability of the system (\ref{equ3.1})--(\ref{equ3.3}), with
solution $\vtwo, \vone, \vnull$, and we again express the
conditionally optimal strategies in terms of the solution to
(\ref{equ3.1})--(\ref{equ3.3}). As mentioned above, this is stronger than only a
verification result.

In Section~\ref{sec4}, we provide equivalent alternative versions for our
BSDEs which are more convenient to work with in some examples with
jumps. This also allows us to discuss in more detail the connections to
the existing literature. Finally, Section~\ref{sec5} illustrates the use
of our results and gives further links to the literature by a number of
simple examples.
%

%
%
\section{Problem formulation and general results}\label{sec1}

We start with a finite time horizon $T\in(0, \infty)$ and a filtered
probability space $(\Omega, \F, \FF, P)$ with the
filtration $\FF= (\F_t)_{0\le t\le T}$ satisfying the usual conditions
of right-continuity and \mbox{$P$-completeness}.
Let $S = (S_t)_{0\le t\le T}$ be an $\RR^d$-valued RCLL semimartingale,
and denote by $\Theta= \Theta_S$ the
space of all predictable $S$-integrable processes $\vartheta$,
$\vartheta\in
L(S)$ for short, such that the stochastic
integral process $\vartheta\SIdot S = \int\vartheta\dd S$ is in the space
$\SS^2(P)$ of semimartingales. Our basic references for terminology and
results from stochastic calculus are \citet{DelMey82} and
\citet{JacShi03}.

For $x\in\RR$ and $H\in L^2(\F_T,P)$, the problem of
\textit{mean-variance hedging} (\textit{MVH}) is to
%
\begin{equation}\label{equ1.1}
\mbox{minimize }E[ (H-x-\vartheta\SIdot S_T)^2 ]\qquad
\mbox{over all
$\vartheta\in\Theta$.}
\end{equation}
The interpretation is that $S$ models the (discounted) prices of $d$
risky assets in a financial market containing also a riskless bank
account with (discounted) price~1. An integrand $\vartheta$ together with
$x\in\RR$ then describes a self-financing dynamic trading strategy with
initial wealth $x$, and $H$ stands for the (discounted) payoff at time
$T$ of some financial instrument. By using $(x,\vartheta)$, we
generate up
to time $T$ via trading a wealth of
$
x + \int_0^T \vartheta_r\dd S_r = x + \vartheta\SIdot S_T
$, and we want to choose $\vartheta$ in such a way that we are close, in
the $L^2(P)$-sense, to the payoff $H$. We embed this into a
\textit{stochastic control problem} and define for $\psi\in\Theta$ and
$t\in[0,T]$
\begin{eqnarray*}
\VH_t(x,\psi)
:\!&=&
\essinf_{\vartheta\in\Theta_{t,T}(\psi)} E[ (H-x-\vartheta
\SIdot S_T)^2
| \F_t ]
\\
&=&
\essinf_{\vartheta\in\Theta_{t,T}(\psi)} E\biggl[ \biggl( H-x-\int
_0^t \psi_r
\dd S_r - \int_t^T \vartheta_r\dd S_r \biggr)^2 \bigg| \F_t
\biggr],
\end{eqnarray*}
where
$
\Theta_{t,T}(\psi):= \{ \vartheta\in\Theta| \vartheta
=\psi$ on $\lsb0, t\rsb\}$.
Our goal is to study the \textit{dynamic value family}
%
\begin{eqnarray}\label{equ1.2}
\VH_t(x):\!&=& \VH_t(x,0)\nonumber\\[-8pt]\\[-8pt]
&=&
\essinf_{\vartheta\in\Theta} E\biggl[ \biggl( H - x - \int_t^T
\vartheta_r\dd S_r
\biggr)^2 \bigg| \F_t \biggr],\qquad
t\in[0,T],
\nonumber
\end{eqnarray}
in order to describe the optimal strategy for the MVH problem (\ref{equ1.1}).
Observe that with this notation, we have the identity
\[
\VH_u\biggl( x+\int_t^u \psi_r\dd S_r \biggr)
=
\VH_u\bigl( x, \psi I_{\rsb t,T\rsb} \bigr)
=
\VH_u\bigl( x, \psi I_{\rsb t,u\rsb} \bigr)
\]
for $u\ge t$. Because the family of random variables
\[
\Gamma_t(\vartheta)
:=
E\biggl[ \biggl( H - x - \int_t^T \vartheta_r\dd S_r \biggr)^2 \bigg|
\F_t \biggr]
\]
for $\vartheta\in\Theta$ is closed under taking maxima and minima,
we have
the classical
\textit{martingale optimality principle} in the following form; see, for
instance, \citet{ElK81} for the general theory, or \citet{ManTev03N1} for a formulation
closer to the present one.
%
\begin{proposition}\label{Prop1.1}
Fix $H\in L^2(\F_T,P)$. For every $x\in
\RR$ and $t\in[0,T]$, we have:

\begin{longlist}[(2)]
\item[(1)]
The process $ ( \VH_u( x + \int
_t^u\vartheta_r\dd S_r ))_{t\le u\le T}$ is a
$P$-submartingale for every \mbox{$\vartheta\in\Theta$}.

\item[(2)] A strategy $\thetast= \thetast(x,H)\in\Theta
_{t,T}(0)$ is
optimal for (\ref{equ1.2}) (i.e., attains the essential
infimum there) if and only if $ ( \VH_u( x + \int
_t^u\thetast_r\dd S_r ))_{t\le u\le T}$ is a \mbox{$P$-martingale}.

\item[(3)] If the strategy $\vartheta^*=\thetasn(x,H)$ solves (\ref{equ1.1}),
then $\vartheta^*
I_{\rsb t,T\rsb}$ is optimal for
$\VH_t(x+\vartheta^*\SIdot S_t) = \VH_t(x, \vartheta^*)$.
\end{longlist}
\end{proposition}

For the special case $H\equiv0$, the fact that $\Theta$ is a cone
immediately gives
%
\begin{equation}\label{equ1.3}
\Vn_t(x)
=
\essinf_{\vartheta\in\Theta} E\biggl[ \biggl( x + \int_t^T \vartheta
_r\dd S_r
\biggr)^2 \bigg| \F_t \biggr]
=
x^2 \Vn_t(1).
\end{equation}
This holds for any random variable $x\in L^2(\F_t,P)$. So Proposition~\ref{Prop1.1}
almost directly gives:
%
\begin{corollary}\label{Cor1.2}
For every $t\in[0,T]$, we have:

\begin{longlist}
\item[(1)] The process $ ( ( 1 + \int
_t^u\vartheta_r\dd S_r
)^2 \Vn_u(1) )_{t\le u\le T}$ is a
$P$-submartingale for every $\vartheta\in\Theta$.

\item[(2)] A strategy $\thetast= \thetast(1,0)\in\Theta
_{t,T}(0)$ is
optimal for $\Vn_t(1)$ in (\ref{equ1.3}) if and only if the process
$ ( ( 1 + \int_t^u\thetast_r\dd S_r )^2 \Vn_u(1)
)_{t\le u\le T}$ is a $P$-martingale.\vadjust{\goodbreak}

\item[(3)] If $\vartheta^*=\thetasn(1,0)$ solves (\ref{equ1.1}) for $x=1$ and
$H\equiv0$, then
%
\begin{equation}\label{equ1.4}
\int_t^T \vartheta^*_r \dd S_r = 0,
\qquad\mbox{$P$-a.s. on the set }\{ 1 + \vartheta^*\SIdot S_t = 0 \}.
\end{equation}
\end{longlist}
\end{corollary}
\begin{pf}
Since (1) and (2) are special cases of Proposition~\ref{Prop1.1}, we
only need to prove (3). Fix $t\in[0,T]$, set $D_t
:= \{ 1 + \vartheta^*\SIdot S_t = 0 \} \in\F_t$ and define $\varphi:=
I_{D_t^c} \vartheta^* I_{\rsb t,T\rsb}$. By part (3) of
Proposition~\ref{Prop1.1} with $x=1, H\equiv0$, the strategy $\vartheta^*
I_{\rsb t,T\rsb}$ is optimal for
$\Vn_t(1+\vartheta^*\SIdot S_t)$ so that
\begin{eqnarray*}
&&
I_{D_t} E\biggl[ \biggl( 1+\vartheta^*\SIdot S_t + \int_t^T \vartheta
^*_r\dd S_r
\biggr)^2 \bigg| \F_t \biggr]\\
&&\qquad
\le
I_{D_t} E\biggl[ \biggl( 1+\vartheta^*\SIdot S_t + \int_t^T \varphi
_r\dd S_r
\biggr)^2 \bigg| \F_t \biggr]
=0
\end{eqnarray*}
by the definitions of $\varphi$ and $D_t$. This yields
\[
0
= I_{D_t} \biggl( 1 + \vartheta^*\SIdot S_t + \int_t^T \vartheta^*_r\dd S_r \biggr)
= I_{D_t} \int_t^T \vartheta^*_r\dd S_r\qquad\mbox{$P$-a.s.}
\]
again by the definition of $D_t$, and so we get (\ref{equ1.4}).
\end{pf}

As in Proposition A.2 of \citet{ManTev03N1}
or Theorem 2.28 of \citet{ElK81}, we also obtain:
%
\begin{proposition}\label{Prop1.3}
Fix $H\in L^2(\F_T,P)$. For every $x\in
\RR$, $t\in[0,T]$ and \mbox{$\psi\in\Theta$}, there exists an
RCLL version of
the $P$-submartingale
\[
\biggl( \VH_u\biggl( x +\int_t^u \psi_r\dd S_r \biggr)
\biggr)_{t\le u\le T}.
\]
Moreover, for each $x\in\RR$, the family $\{\VH_t(x) | t\in
[0,T] \}
$ of random variables can be aggregated into an RCLL process, which we
again call
\[
\VH(x) =( \VH_u(x))_{0\le u\le T}.
\]
\end{proposition}

In the sequel, we always choose and work with the RCLL versions from
Proposition~\ref{Prop1.3}.

For easier discussion of the next result, we introduce some more
terminology. We denote by $\PPes$
the (a priori possibly empty) set of all probability measures $Q$
equivalent to $P$ on $\F_T$ such that $S$ is a
$Q$-$\sigma$-martingale and ${dQ\over dP}\in L^2(P)$. Assuming that
$\PPes$ is nonempty is one way of imposing \textit{absence of arbitrage}
for our financial market and also fits naturally with the fact that our
basic problem is cast in quadratic terms. The density process of $Q$
with respect to $P$ is denoted by $\ZQ= (\ZQ_t)_{0\le t\le T}$, and we
say that $Q\in\PPes$ satisfies the reverse H\"older inequality
$R_2(P)$ if there is a constant $C$ with\vadjust{\goodbreak}
$
E_P[ (\ZQ_T)^2 | \F_\tau]
\le
C (\ZQ_\tau)^2
$
$P$-a.s. for all stopping times $\tau\le T$. It is well known that if there
is some $Q\in\PPes$ satisfying $R_2(P)$, then $G_T(\Theta) =\{
\vartheta
\SIdot S_T | \vartheta\in\Theta\}$ as well as $L^2(\F_t,P) + G_T(
\Theta_{t,T}(0))$ for each $t$ are closed in $L^2(P)$ so that both
(\ref{equ1.1}) and (\ref{equ1.2}) for each $t$ have a solution; see Theorem
5.2 of \citet{ChoKraStr98}. Moreover, for any $Q\in
\PPes$ and any $\vartheta\in\Theta$,
the product of $\ZQ$ and
$\vartheta\SIdot S$ is a $P$-$\sigma$-martingale with $P$-integrable
supremum; so $\vartheta\SIdot S$ is a true
$Q$-martingale, and $\vartheta\SIdot S_T=0$ a.s. implies that
$\vartheta=0$
in $L(S)$. This is used later several times to argue that a
self-financing strategy is uniquely determined by its wealth process
(i.e., stochastic integral).

Our main result in this section now provides the basic structure of the
process $\VH(x)$ and of the optimal
strategies for (\ref{equ1.2}).
%
\begin{theorem}\label{Theorem1.4}
Fix $H\in L^2(\F_T,P)$. Suppose that for
each $t\in[0,T]$, (\ref{equ1.2}) has a solution
$\thetast= \thetast(x,H)$ for every $x\in\RR$. Suppose also that for
any $\vartheta\in\Theta$, $\vartheta\SIdot S_T=0$ a.s. implies that
$\vartheta
=0$ in $L(S)$. Then each $\thetast(x,H)$ is of the
affine form
%
\begin{equation}\label{equ1.5}
\thetast(x,H) = \vartheta^{0,t} + x \vartheta^{1,t}
\qquad\mbox{for some }\vartheta^{0,t}, \vartheta^{1,t}\in\Theta_{t,T}(0),
\end{equation}
and each $\VH_t(x)$ has the quadratic form
%
\begin{equation}\label{equ1.6}
\VH_t(x) = \vnull_t - 2\vone_t x + \vtwo_t x^2
\end{equation}
for RCLL processes $\vnull, \vone, \vtwo$ not depending on $x$.
Moreover, $\vartheta^{1,t} = \thetast(1,0)$ is the solution of (\ref{equ1.3}),
and the quadratic coefficient $\vtwo_t$ equals $\Vn_t(1)$ from
(\ref{equ1.3}) and does not depend on $H$.
\end{theorem}
\begin{pf}
Fix\vspace*{1pt} $t\in[0,T]$. Denote by $G_{t,T} = G_T(
\Theta_{t,T}(0)) = \{ \int_t^T \vartheta_r\dd S_r |
\vartheta\in\Theta\}$ the space of all stochastic integrals on $\rsb
t,T\rsb$ of $\vartheta\in\Theta$ and by $\bar{G}_{t,T}$ its closure in
$L^2(P)$. Since the problems (\ref{equ1.2}) with payoff $H$ for $x=1$
and $x=0$ have solutions (which are given by projections), so does
problem (\ref{equ1.2}) for $x=1$ and payoff \mbox{$H^\prime\equiv0$} by
taking differences, and the latter problem is identical to
(\ref{equ1.2}) for $x=0, H^\prime\equiv-1$ so that $\vartheta^*(0,-1) =
\vartheta^*(1,0)$. Both here and in the next argument, we exploit our
assumption that a self-financing strategy\vspace*{1pt} is uniquely
determined by its wealth process. If $\Pi$ is the projection in
$L^2(P)$ on $\bar{G}_{t,T}$, then clearly
\begin{eqnarray*}
\thetast(x,H)\SIdot S_T
&=&
\Pi(H-x)
=
\Pi(H) + x \Pi(-1)\\
&=&
\thetast(0,H)\SIdot S_T + x \thetast(0,-1)\SIdot S_T,
\end{eqnarray*}
and so (\ref{equ1.5}) follows with $\vartheta^{0,t} = \thetast(0,H)$ and
$\vartheta^{1,t} = \thetast(0,-1) = \thetast(1,0)$. This gives
\begin{eqnarray*}
\VH_t(x)
&=&
E\biggl[ \biggl( H-x-\int_t^T \thetast_r(x,H) \dd S_r \biggr)^2 \bigg|
\F
_t \biggr]\\
&=&
E\biggl[ \biggl( H - \int_t^T\vartheta^{0,t}_r\dd S_r - x\biggl( 1+\int_t^T
\vartheta^{1,t}_r\dd S_r \biggr) \biggr)^2 \bigg|
\F_t \biggr],
\end{eqnarray*}
and hence we directly obtain the expression (\ref{equ1.6}) with
\begin{eqnarray*}
\vnull_t
&=&
E\biggl[ \biggl( H-\int_t^T \thetast_r(0,\allowbreak H) \dd S_r
\biggr)^2 \bigg| \F_t
\biggr],
\\
\vone_t
&=&
E\biggl[ \biggl( H-\int_t^T \thetast_r(0, H) \dd S_r \biggr) \biggl(
1+\int_t^T
\thetast_r(1,0) \dd S_r \biggr) \bigg| \F_t \biggr]
\end{eqnarray*}
and
%
\begin{equation}\label{equ1.7}
\vtwo_t
=
E\biggl[ \biggl( 1+\int_t^T \thetast_r(1,0) \dd S_r \biggr)^2 \bigg|
\F_t
\biggr] = \Vn_t(1).
\end{equation}
Since the families $\{\VH_t(x) | t\in[0,T] \}$ aggregate into an
RCLL process, the same holds for the families $\vnull, \vone, \vtwo$
from (\ref{equ1.6}). The last assertion is clear from the above proof.
\end{pf}
\begin{Rems*}
(1) As mentioned above, one sufficient condition for all assumptions
of Theorem~\ref{Theorem1.4} is the existence of some $Q\in\PPes$ satisfying
the reverse H\"older inequality $R_2(P)$; see \citet{ChoKraStr98}.

(2) The particular choice of $\Theta= \Theta_S$ for the
space of
integrands is convenient and also exploited later, but not crucially
important for the conclusion of Theorem~\ref{Theorem1.4} to hold. All we need
is that there exist for all $t$ solutions $\thetast(x,H)$ for all $x$,
that the martingale optimality principle from Proposition~\ref{Prop1.1}
holds, and that $\Theta$ [or $G_T(\Theta)$, which must be a subset of
$L^2(P)$] is a linear space. Of course, existence of solutions for all
$x$ and all $H$ is equivalent to closedness of $G_T(\Theta)$ in
$L^2(P)$; and the key point for the martingale optimality principle is
closedness under bifurcation of $\Theta$.

(3) We emphasize that Theorem~\ref{Theorem1.4} is a bit of a folklore
result in the literature on mean-variance hedging, and we do not claim
any great originality here. Variants in different levels of generality
can be found in \citet{Gug03}, \citet{ManTev03N1},
\citet{BobSch04}, \citet{autokey5}, to name but a few.
However, we think that it is useful to have a presentation which is as
general, and yet as simple, as possible.
\end{Rems*}

Our goal in the sequel is to study the dynamics of the coefficient
processes $\vnull, \vone, \vtwo$ and use them to express the optimal
strategies $\thetast(x,H)$. Let us first simplify things a little.
Because $\thetast(1,0)$ is the solution (minimizer) of (\ref{equ1.3}), the
first order condition for that quadratic optimization
problem implies that
$
E[ \int_t^T \vartheta_r\dd S_r ( 1+\int_t^T \thetast
_r(1,0) \dd S_r ) | \F_t ] = 0
$
$P$-a.s. for each $t\in[0,T]$ and $\vartheta\in\Theta$. We note for
later use
that this allows us to write
%
\begin{equation}\label{equ1.8}
\vone_t = E\biggl[ H \biggl( 1+\int_t^T \thetast_r(1,0) \dd S_r \biggr)
\bigg| \F_t \biggr].
\end{equation}

Also for later use, we give some additional results for the
coefficients $\vnull, \vone$, $\vtwo$.
%
\begin{lemma}\label{Lem1.5}
Under the assumptions of Theorem~\ref{Theorem1.4}, we have:

\begin{longlist}[(2)]
\item[(1)] $\vtwo$ is a $P$-submartingale with $0\le\vtwo\le1$.

\item[(2)] $\vnull$ is a $P$-submartingale with $0\le\vnull_t\le
E[ H^2
| \F_t ]$, $0\le
t\le T$, hence of class (D).

\item[(3)] $\vone$ is a $P$-special semimartingale with $|\vone
|^2$ of
class (D). Therefore $\vone$ is in $\SS^2_\loc(P)$ and for its
canonical decomposition $\vone= \vone_0 + \mone+ \aone$, we
have $\mone\in\M^2_{0,\loc}(P)$.
\end{longlist}
\end{lemma}
\begin{pf}
(1) By Theorem~\ref{Theorem1.4} and (\ref{equ1.7}), we have $\vtwo= \Vn
(1)$, and this is a $P$-submartingale by part (1) of Corollary~\ref{Cor1.2} (for $\vartheta\equiv0$). Because $\vartheta\equiv0$ is in
$\Theta
$, we get $0\le\Vn(1)\le1$ directly from
(\ref{equ1.3}).

(2) Theorem~\ref{Theorem1.4} gives $\vnull= \VH(0)$, and this is a
$P$-submartingale by part~(1) of Proposition~\ref{Prop1.1}
(for $x=0, \vartheta\equiv0$) and nonnegative by the definition in
(\ref{equ1.2}). Since $\vartheta\equiv0$ is in $\Theta$, (\ref{equ1.2})
also gives $\VH_t(0) \le E[ H^2 | \F_t ]$ for all~$t$.

(3) By part (1) of Proposition~\ref{Prop1.1}, $\VH(x)$ is a
$P$-submartingale, hence a \mbox{$P$-special} semimartingale, and so
are\vspace*{1pt} $\vtwo$ and $\vnull$ by (1) and (2). Because $\VH(x) = \vnull-
2\vone x + \vtwo x^2$ by Theorem~\ref{Theorem1.4},
also $\vone$ is then a $P$-special semimartingale. Moreover,\vspace*{1pt} $\VH
(x)\ge
0$ for all $x$ due to (\ref{equ1.2}) implies that
\[
\bigl|\vone_t\bigr|^2 \le\vtwo_t\vnull_t \le\vnull_t \le E[ H^2 | \F_t ],
\qquad
0\le t\le T,
\]
by (1) and (2) so that $|\vone|^2$ is of class (D). The
rest of part (3) is then clear.
\end{pf}

%
%
\section{\texorpdfstring{Pure investment: The special case $x=1$, $H\equiv0$}
{Pure investment: The special case x=1, H equivalent 0}}\label{sec2}

In this section, we give a description of (the RCLL version of) the
value process
%
\begin{equation}\label{equ2.1}
\Vn_t(1)
=
\essinf_{\vartheta\in\Theta} E\biggl[ \biggl( 1+\int_t^T \vartheta
_r\dd S_r
\biggr)^2 \bigg| \F_t \biggr],
\qquad0\le t\le T,
\end{equation}
of the problem (\ref{equ1.3}). Since this is by (\ref{equ1.7}) and Theorem
\ref{Theorem1.4} the \textit{quadratic} coefficient in the
representation (\ref{equ1.6}), we use in this section the shorter notation
%
%
\[
\q_t:= \Vn_t(1) = \vtwo_t,
\qquad 0\le t\le T.
\]
We also remark that $\q$ coincides with the \textit{opportunity process}
from \citet{CerKal07}, although the latter is
defined there with a different space $\Theta$ of integrands $\vartheta$
for $S$.

Let us first prove strict positivity of $\q$, as well as of $\q_-$.
%
\begin{lemma}\label{Lem2.1}
Suppose $\PPes\ne\varnothing$. Then $\q$ and $\q
_-$ are both strictly positive, in the sense that
$P[\q_t>0$ and $\q_{t-}>0$ for $0\le t\le T] = 1$.
If there is some $Q\in\PPes$ satisfying the reverse H\"older inequality
$R_2(P)$, we even have
$\q\ge\delta>0$ $P$-a.s. for some constant $\delta$.
\end{lemma}
\begin{pf}
For $Q\in\PPes$ with density process $Z = \ZQ= \ZQP$, define as
in \citet{GouLauPha98} a new probability
$R\approx P$ by
$
{dR\over dP}:= {Z_T^2 \over E[Z_T^2]}
$.
Then the Bayes rule gives
%
\begin{eqnarray}\label{equ2.2}
\ZRP_t
&:=&
{dR\over dP}\bigg|_{\F_t} = {E[Z_T^2 | \F_t] \over E[Z_T^2]},
\\
\label{equ2.3}
\ZRQ_t
&:=&
{dR\over dQ}\bigg|_{\F_t} = {E_Q[Z_T | \F_t] \over E[Z_T^2]} =
{1\over Z_t} \ZRP_t.
\end{eqnarray}
Using the Bayes rule and (\ref{equ2.2}), Jensen's inequality, again the
Bayes rule and (\ref{equ2.3}) yields
\begin{eqnarray*}
&&
E\biggl[ \biggl( 1+\int_t^T \vartheta_r\dd S_r \biggr)^2 \bigg| \F
_t \biggr]
%
\\
&&\qquad =
\ZRP_t E[Z_T^2] E_R\biggl[ ( Z_T^2)^{-1} \biggl( 1+\int_t^T \vartheta
_r\dd S_r
\biggr)^2 \bigg| \F_t \biggr]
\\
&&\qquad \ge
\ZRP_t E[Z_T^2] \biggl( E_R\biggl[ (Z_T)^{-1} \biggl( 1+\int_t^T
\vartheta_r\dd S_r \biggr) \bigg| \F_t \biggr]
\biggr)^2
\\
&&\qquad =
\ZRP_t E[Z_T^2] \biggl( ( \ZRQ_t)^{-1} E_Q\biggl[ ( E[Z_T^2])^{-1} \biggl(
1+\int_t^T \vartheta_r\dd S_r
\biggr) \bigg| \F_t \biggr] \biggr)^2.
\end{eqnarray*}
But as already noted before Theorem~\ref{Theorem1.4}, $\int\vartheta\dd S$ is a
$Q$-martingale whenever $Q\in\PPes$ and
$\vartheta\in\Theta$. So we get by using (\ref{equ2.3}) and (\ref{equ2.2}) that
%
\begin{equation}\label{equ2.4}
E\biggl[ \biggl( 1+\int_t^T \vartheta_r\dd S_r \biggr)^2 \bigg| \F
_t \biggr]
\ge
{\ZRP_t E[Z_T^2] \over( \ZRQ_t E[Z_T^2])^2}
=
{Z_t^2 \over E[Z_T^2 | \F_t]},
\end{equation}
and the first assertion follows since $\inf_{0\le t\le T} Z_t > 0$
$P$-a.s. by the minimum principle for supermartingales and $\sup_{0\le
t\le T} E[Z_T^2 | \F_t] < \infty$ $P$-a.s. by the\vspace*{1pt}
martingale maximal inequality. If $Q$ satisfies $R_2(P)$ with constant
$C$, we can take $\delta= 1/C$ for the second claim.
\end{pf}
\begin{Rem*}
Strict positivity of the opportunity process and its left
limits (hence of $\q$ and $\q_-$) is also
proved in Lemma 3.10 of \citet{CerKal07}. However, the
above short proof seems to us
more transparent.
\end{Rem*}

The optimization problem in (\ref{equ2.1}) has a (well-known) dual formulation
as follows. Extending $\PPes$ a little, we denote by $\PPss$ the set of
all \textit{signed} measures $Q\ll P$ on $\F _T$ with $Q[\Omega] = 1$
and such that the product of $S$ and the density process $\ZQ$ of $Q$
with respect to $P$ is a $P$-$\sigma$-martingale. We call\vspace*{-1pt} $\Qt\in\PPss$
\textit{variance-optimal} if $ \Vert{d\Qt\over dP}\Vert_{L^2(P)} \le
\Vert{dQ\over dP}\Vert_{L^2(P)} $ for all $Q\in\PPss$, and we say that
the \textit{variance-optimal martingale measure} (\textit{VOMM}) exists
if $\Qt\in\PPes$ is variance-optimal. (In particular, $\Qt$ is then by
definition equivalent to $P$.) If $S$ is continuous, Theorem 1.3 of
\citet{DelSch96} shows that $\PPes\ne\varnothing$ is
sufficient for the VOMM to exist; but if $S$ can have jumps, the
situation is more complicated.

The dynamic problem of finding the VOMM has the value process
\[
\Vt_t:= \essinf_{Q\in\PPes} E[ ( \ZQ_T / \ZQ_t
)^2
| \F_t ],\qquad
0\le t\le T.
\]
Then we have the following direct connection to $\Vn(1)$ and (\ref{equ2.1}).
%
\begin{proposition}\label{Prop2.2}
Suppose $S\in\SS^2_\loc(P)$ and that
the VOMM exists. Then $\Vt= 1/\Vn(1)$.
\end{proposition}
\begin{pf}
We know from (\ref{equ2.4}) in the proof of Lemma~\ref{Lem2.1} that
for $\vartheta\in\Theta$ and $Q\in\PPes$,
\[
E\biggl[ \biggl( 1+\int_t^T \vartheta_r\dd S_r \biggr)^2 \bigg| \F
_t \biggr]
\ge
1 / E[ ( \ZQ_T / \ZQ_t )^2 | \F_t ],
\qquad 0\le t\le T.
\]
Taking the ess inf over $\vartheta\in\Theta$ and the ess sup over
$Q\in
\PPes$ implies that $\Vn(1) \ge1/\Vt$.
Conversely, 
since $\Vn_T(1) = 1$, the martingale optimality principle
in Corollary~\ref{Cor1.2} gives
%
\begin{equation}\label{equ2.5}\qquad
\biggl( 1 + \int_0^t \vartheta_r\dd S_r \biggr)^2 \Vn_t(1)
\le
E\biggl[ \biggl( 1+\int_0^T \vartheta_r\dd S_r \biggr)^2 \bigg| \F
_t \biggr],\qquad
0\le t\le T,
\end{equation}
for every $\vartheta\in\Theta=\Theta_S$. But if we define, as in
\citet{GouLauPha98},
\begin{eqnarray*}
\Theta_{\mathrm{GLP}}
&:=&
\{ \vartheta\in L(S) | \vartheta\SIdot S_T\in L^2(P)\mbox{ and }
\ZQ(\vartheta\SIdot S)\\
&&\hspace*{20pt}\mbox{is a
$P$-martingale for all $Q\in\PPss$}\},
\end{eqnarray*}
%
%
then $G_T(\Theta_{\mathrm{GLP}}):= \{\vartheta \cdot S_T | \vartheta
\in \Theta_{\mathrm{GLP}}\}$ is by Corollary 2.9 of
\citet{CerKal07} the closure of $G_T(\Theta_S)$ in $L^2(P)$,
and this allows us to extend
(\ref{equ2.5}) to every $\vartheta\in\Theta_{\mathrm{GLP}}$. Indeed,
for a sequence $(\vartheta ^n)$ in $\Theta_S$ with $G_T(\vartheta^n)
\to G_T(\vartheta)$ in $L^2(P)$, the right-hand side of (\ref{equ2.5})
for $\vartheta^n$ converges in $L^1(P)$ to the right-hand\vspace*{1pt} side of
(\ref{equ2.5}) for $\vartheta$, and because we have $ \ZQ_t ( 1 +
\int_0^t \vartheta_r\dd S_r ) = E[ \ZQ_T ( 1 + \int_0^T \vartheta_r\dd
S_r ) | \F_t ] $ for
$\vartheta\in\Theta_{\mathrm{GLP}}\supseteq\Theta_S$ and $Q\in \PPes$,
the left-hand side of (\ref{equ2.5}) for $\vartheta^n$ converges in
probability to the left-hand side of (\ref{equ2.5}) for $\vartheta$. We
remark that the use of Corollary 2.9 in \citet{CerKal07} exploits
that $S\in\SS ^2_\loc(P)$.\vspace*{2pt}

By assumption, the VOMM $\Qt$ exists. A slight modification of the
proof of Lemma 2.2 in \citet{DelSch96} [since $S$ is in
$\SS^2_\loc(P)$ instead of locally bounded] yields $ Z^{\Qt}_T = c +
\int_0^T \thetat_r\dd S_r $ for some $c>0$ and
$\thetat\in\Theta_{\mathrm{GLP}}$ and thus $ E_{\Qt}[ Z^{\Qt}_T | \F_t]
= c + \int_0^t \thetat_r\dd S_r, $ $0\le t\le T$. Applying\vspace*{1pt}
(\ref{equ2.5}) with $\vartheta:= \thetat/c$ and using the Bayes rule
therefore gives
\[
( Z^{\Qt}_t)^2 E[ ( Z^{\Qt}_T )^2 | \F_t
]
\ge
( Z^{\Qt}_t)^2 ( E_{\Qt} [ Z^{\Qt}_T | \F_t
]
)^2 \Vn_t(1)
=
( E[ ( Z^{\Qt}_T )^2 | \F_t ]
)^2 \Vn_t(1)
\]
and hence
\[
1/ \Vn_t(1) \ge E[ ( Z^{\Qt}_T / Z^{\Qt}_t )^2
| \F
_t ] \ge\Vt_t,
\qquad 0\le t\le T.
\]
This completes the proof.
\end{pf}
\begin{Rem*}
For experts on mean-variance hedging, Proposition~\ref{Prop2.2}
is also a kind of folklore result. For the case where the filtration
is continuous, it can, for instance, be found in Proposition 4.2 of
\citet{ManTev03N1} (with the remark that it extends to general
$\FF$ if $S$ is
continuous). But we do not know a reference for the level of generality
given here.
\end{Rem*}

Henceforth, we often use the following simple fact:
%
\begin{equation}\label{equ2.6}
\begin{tabular}{p{243pt}}
If $B,C$ are of locally integrable variation and $B\ll
C$,
then
also $\dpP{B} \ll\dpP{C}$.
\end{tabular}
\end{equation}
In (\ref{equ2.6}), the (right) superscript ${}^{\p}$ denotes the
compensator or dual predictable projection. This should not be confused
with the predictable projection of a process $Y$ which is denoted by
$\pY$, with a left superscript. The most frequent application of
(\ref{equ2.6}) will be for $C = [M]$, where $\dpP{C} = \dpP{[M]} = \MM$
when $M$ is a locally square-integrable local
martingale.

\textit{In the sequel}, \textit{we focus on the case $d=1$ so that $S$ is
one-dimensional.} One can obtain analogous results for $d>1$ (and we
shall comment on this later), but the arguments and formulations look
more technical without providing extra insight. When
$S\in\SS^2_\loc(P)$ so that $S$ is in particular a $P$-special
semimartingale, we write $S = S_0+M+A$ for its $P$-canonical
decomposition and note that $M\in\M ^2_{0,\loc}(P)$ and $A$ is
predictable and of locally square-integrable (or even locally bounded)
variation. If we also have $\PPes\ne\varnothing$, then it is well known
that $S$ satisfies the so-called structure condition, that is, that $S$
has the form
%
\begin{equation}\label{equ2.7}
S = S_0 + M + A = S_0 + M + \int\lambda\dd\MM
\end{equation}
with $M\in\M^2_{0,\loc}(P)$ and $\lambda\in L^2_\loc(M)$; see
Theorem 1
of \citet{Sch95}. This implies that
\[
[A]
=
\biggl[ \int\lambda\dd\MM\biggr]
=
\sum(\lambda_s \Delta\MM_s)^2
=
(\lambda^2 \Delta\MM)\SIdot\MM
\ll
\MM.
\]
Because $A$ is predictable, $[M,A]$ is a local $P$-martingale by
Yoeurp's lemma so that
%
\begin{equation}\label{equ2.8}
\dpP{[S]} = \dpP{([M] + [A])} = (1+\lambda^2\Delta\MM)\SIdot\MM.
\end{equation}

Now suppose that $S\in\SS^2_\loc(P)$ and $\PPes\ne\varnothing$. To
describe the process $\q=\Vn(1)$ by a BSDE, we first introduce an
auxiliary operation. Suppose $Y$ is a $P$-special semimartingale with
canonical decomposition $Y = Y_0 + \NY+ B^Y$. Then
$
[ Y, [S] ]
=
[ \NY, [S] ] + \Delta B^Y\SIdot[S]
$,
and if $[ \NY, [S] ]$ is of locally $P$-integrable variation, we have
by (\ref{equ2.8}) and (\ref{equ2.6}) that
%
\begin{equation}\label{equ2.9}
\dpP{\bigl[ Y, [S] \bigr]} = \dpP{\bigl[ \NY, [S] \bigr]} + \Delta B^Y\SIdot\dpP{[S]} \ll
\MM.
\end{equation}
Note also that the predictable stopping theorem gives $\Delta B^Y =
{}^{\p}\Delta Y = \pY- Y_-$ so that
%
\begin{equation}\label{equ2.10}
Y_- + \Delta B^Y = \pY.
\end{equation}
The auxiliary quantity we need is the predictable
Radon--Nikod\'ym derivative
%
\begin{equation}\label{equ2.11}
\func_t(Y):= {d \dpP{[\NY, [S] ]}_t \over d\MM_t},
\qquad 0\le t\le T.
\end{equation}
Finally, we introduce the notation
%
\begin{equation}\label{equ2.12}
{\mathcal N}(Y):= \pY(1+\lambda^2\Delta\MM) + \func(Y).
\end{equation}
The condition that $[ \NY, [S] ]$ is in $\A_\loc(P)$ (and hence
has a
compensator) is, for instance, satisfied if $Y$ is bounded, hence in
particular for $Y=\q$.
\begin{Rem*}
In the context of the equations we study, the operation ${\mathcal
N}(Y)$ in (\ref{equ2.12}) can sometimes be simplified. If $S$ is continuous,
then so are $[S]$ and $\MM$, due to (\ref{equ2.7}); so $\func(Y)$ and
$\Delta\MM$ then both vanish and (\ref{equ2.12}) reduces to the expression ${\mathcal N}(Y) =
\pY= Y_- + \Delta B^Y$. Looking ahead at (\ref{equ2.18}), however, we see that
we are interested in the case where $B^Y\ll\MM$, and so we then also
get $\Delta B^Y = 0$ and hence $\NY= Y_-$. Finally, if even the
filtration $\FF$ is continuous, then $L$ in (\ref{equ2.18}) is continuous; so is
then $Y$, and we end up with ${\mathcal N}(Y) = Y$.
\end{Rem*}

Our next result shows that ${\mathcal N}(\q) = \N(\vtwo)$ is always strictly
positive. This is important since we later need to divide by ${\mathcal
N}(\q)$.
%
\begin{lemma}\label{Lem2.3}
Suppose $\PPes\ne\varnothing$ and $S\in\SS
^2_\loc(P)$. If $\q\ge\delta>0$
for some constant $\delta$, then
%
\begin{equation}\label{equ2.13}\qquad
{\mathcal N}(\q) = \pq(1+\lambda^2\Delta\MM) + \func(\q) \ge\delta,
\qquad\mbox{$P\otimes\MM$-a.e. on $\lsb0, T\rsb$.}
\end{equation}
In general, we still have
%
\begin{equation}\label{equ2.14}
{\mathcal N}(\q) > 0,
\qquad\mbox{$P\otimes\MM$-a.e. on $\lsb0, T\rsb$.}
\end{equation}
Moreover, ${\mathcal N}(\q)$ is locally bounded away from 0 (uniformly in
$t,\omega$).
\end{lemma}
\begin{pf}
If $\q\ge\delta$, then $B:= \q\SIdot[S] - \delta[S]$ is in $\A
^+_\loc(P)$ and hence also
$\dpP{B}\in\A^+_\loc(P)$. But $B\ll[S]$, hence $\dpP{B} \ll\dpP
{[S]} =
(1+\lambda^2\Delta\MM) \SIdot\MM$ by (\ref{equ2.6}) and (\ref{equ2.8}), and so
\begin{eqnarray*}
\dpP{B}
&=&
\dpP{( \q\SIdot[S] )} - \delta(1+\lambda^2\Delta\MM) \SIdot\MM
\\
&=&
\int\biggl( {d \sdpP{( \q\sSIdot[S] )} \over d\MM} -
\delta(1+\lambda
^2\Delta\MM) \biggr) \dd\MM
\in\A^+_\loc(P).
\end{eqnarray*}
Writing $\q= \q_- + \Delta\q$ and $\Delta\q\SIdot[S] = [ \q,
[S] ]$
and using (\ref{equ2.8})--(\ref{equ2.12}) yields
%
\begin{eqnarray}
\label{equ2.15}
\dpP{( \q\SIdot[S] )}
&=&
\q_- \SIdot\dpP{[S]} + \dpP{[ \Nq, [S] ]} + \Delta\Bq
\SIdot
\dpP{[S]}
\nonumber\\
&=&
\bigl( \pq(1+\lambda^2\Delta\MM) + \func(\q) \bigr) \SIdot\MM
\\
&=&
{\mathcal N}(\q) \SIdot\MM.\nonumber
\end{eqnarray}
Thus we obtain $ \dpP{B} = \{ {\mathcal N}(\q) -
\delta(1+\lambda^2\Delta\MM) \} \SIdot \MM \in\A ^+_\loc(P) $, and this
implies (\ref{equ2.13}) since $\lambda^2\Delta\MM\ge0$. In general,
setting $ \tau_n:= \inf\{ t\ge0 | \q_t < {1\over n}\} \wedge T $ (with
$\inf\varnothing= +\infty$) gives $\tau_n\nearrow T$ stationarily
because $\q>0$ by Lemma~\ref{Lem2.1}, and $\q\ge{1\over n}$ on $D_n:=
\lsb0, {\tau_n\lsb}\cup(\Omega\times\{T\})$ since $\q_T=1$. The
argument for (\ref{equ2.13}) now implies that ${\mathcal N}(\q)
\ge{1\over n}$ holds $P\otimes\MM$-a.e. on $D_n$, and (\ref{equ2.14})
follows since $\bigcup _{n\in\NN} D_n = \lsb0,T\rsb$.
For the final assertion, note that the preceding proof shows that
${\mathcal N} (\q)^{\tau_n-} \ge{1\over n}$ so that the nonnegative
process $1/{\mathcal N}(\q)$ is prelocally bounded. Since $1/{\mathcal
N}(\q)$ is like ${\mathcal N}(\q)$ predictable, it is therefore by
\citet{DelMey82}, Remark VIII.11 also locally bounded, and this means
that ${\mathcal N}(\q)$ is locally bounded away from 0.
\end{pf}
\begin{Rem*}
If $d>1$, both $[S]$ and $\MM$ have to be replaced by
matrix-valued processes $( [S^i, S^j])_{i,j=1, \ldots, d}$ and $(
\langle M^i, M^j \rangle)_{i,j=1, \ldots, d}$. We then take a
predictable $B\in\A^+_\loc(P)$ with $\langle M^i, M^j\rangle= \mu
^{ij}\SIdot B \ll B$ and define the matrix-valued predictable process
$\func(\q)$ by
%
\begin{equation}\label{equ2.16}
\func^{ij}_t(\q):= {d \dpP{[\Nq, [S^i, S^j] ]}_t \over dB_t},
\qquad 0\le t\le T.
\end{equation}
Analogously to Lemma~\ref{Lem2.3}, one can then prove that
%
\begin{equation}\label{equ2.17}
{\mathcal N}(\q):= \pq\bigl( \mu+ (\mu\lambda)\tr\mu\lambda
\Delta B
\bigr) + \func(\q) \mbox{ is positive definite $P\otimes B$-a.e.}\hspace*{-30pt}
\end{equation}
\end{Rem*}

Recalling the notation (\ref{equ2.12}), we now consider the \textit{backward equation}
%
\begin{eqnarray}
\label{equ2.18}
Y_t
&=&
Y_0 + \int_0^t {( \psi_s + \lambda_s \pY_s)^2 \over\NNY{s} } \dd\MM_s +
\int_0^t \psi_s\dd M_s + L_t
\nonumber\hspace*{-29pt}\\
&=&
Y_0 + \int_0^t {( \psi_s + \lambda_s \pY_s)^2 \over\pY_s
(1+\lambda
_s^2 \Delta\MM_s) + \func_s(Y) } \dd\MM_s
+ \int_0^t \psi_s\dd M_s + L_t,\hspace*{-29pt}\\
Y_T&=&1.\nonumber\hspace*{-29pt}
\end{eqnarray}
A \textit{solution} of (\ref{equ2.18}) is a triple $(Y, \psi, L)$, where $L$
is a local $P$-martingale which is strongly $P$-orthogonal to $M$,
$\psi
$ is in $L^1_\loc(M)$ and $Y = Y_0 + \NY+ B^Y$ is a \mbox{$P$-special}
semimartingale with $[ \NY, [S] ] \in\A_\loc(P)$. Note that
$\lambda$
and $M$ come from $S$ via (\ref{equ2.7}). With a slight abuse of
terminology, we sometimes call $Y$ instead of the whole triple $(Y,
\psi
, L)$ a solution; any properties then only refer to $Y$.

Denoting the stochastic exponential started at time $t$ of a
semimartingale $X$ by
\[
\tE(X)_u = 1 + \int_t^u \tE(X)_{r-} \dd X_r = \E( X - X^t)_u,
\qquad t\le u\le T,
\]
our first main result is the following description of $\Vn(1) = \q$ via
a BSDE.
%
\begin{theorem}\label{Theorem2.4}
Suppose that $S\in\SS^2_\loc(P)$ and $\PPes
\ne\varnothing$. Then:

\begin{enumerate}[(1)]
\item[(1)] The following two assertions are equivalent:

\begin{enumerate}[(a)]
\item[(a)] For every $t\in[0,T]$, there exists an optimal
strategy $\thetast(1,0)\in\Theta_{t,T}(0)$ for (\ref{equ1.2}) with
$x=1, H\equiv0$.

\item[(b)] There exists a solution $(Y, \psi, L)$ to the BSDE
(\ref{equ2.18}) having $L\in\M^2_{0,\loc}(P)$, $\psi\in L^2_\loc(M)$, $Y$
bounded and strictly positive and such that for every $t\in[0,T]$,
the process $( \tE( - {\psi+ \lambda\spP{Y} \over{\mathcal N}(Y)}
\SIdot S)_u
)_{t\le u\le T}$ is in $\SS^2(P)$.
\end{enumerate}
If \textup{(a)} or \textup{(b)} hold, then the optimal $\thetast(1,0)$ is for every $t$
given by
%
\begin{equation}\label{equ2.19}
\thetast_u(1,0) = - {\psi_u + \lambda_u \pY_u \over\NNY{u}} \tE \biggl( -
{\psi+ \lambda\pY\over{\mathcal N}(Y)} \SIdot S \biggr)_{u-},\qquad
t\le u\le T,
\end{equation}
and $\q= \Vn(1)$ is the unique bounded strictly positive solution of
(\ref{equ2.18}).

\item[(2)] Suppose, in addition, that there is some $Q\in
\PPes$ satisfying
the reverse H\"older inequality
$R_2(P)$. Then $\q= \Vn(1)$ is the unique solution to the BSDE (\ref{equ2.18})
in the class of processes satisfying $c\le Y\le
C$ for positive constants~$c,C$. Moreover, the optimal $\thetast(1,0)$
exist and are given by (\ref{equ2.19}).
\end{enumerate}
\end{theorem}
\begin{pf}
Throughout this proof, we write $\thetast$ for $\thetast(1,0)$
and denote by $m$ a generic local $P$-martingale that can change from
one appearance to the next.\vspace*{8pt}

(1) For part (1) of Theorem~\ref{Theorem2.4}, we start by deriving the
BSDE (\ref{equ2.18}). By part (1) of Lemma~\ref{Lem1.5}, $\q= \vtwo$ is
a $P$-submartingale, hence a $P$-special semimartingale with canonical
decomposition $ \q= \q_0 + \Nq+ \Bq$, and $0\le\q\le1$ implies that
$\q\in\SS^2_\loc(P)$ and $\Nq$ has bounded jumps and is in $\M^2_{0,
\loc}(P)$. The Galtchouk--Kunita--Watanabe decomposition thus allows us
to write
%
\begin{equation}\label{equ2.20}
\q= \q_0 + \varphi\SIdot M + L^{\q} + \Bq
\end{equation}
with $\varphi\in L^2_\loc(M)$ and $L^{\q}\in\M^2_{0,\loc}(P)$ strongly
$P$-orthogonal to $M$. Combining this with (\ref{equ2.7}) and Yoeurp's
lemma then gives
%
\begin{eqnarray}\label{equ2.21}
[\q,S] &=& m + \varphi\SIdot[M] + [A, \Bq] \nonumber\\[-8pt]\\[-8pt]
&=& m + ( \varphi+
\lambda\Delta\Bq)\SIdot\MM.\nonumber
\end{eqnarray}

We now apply It\^o's formula to the process
$
X^\vartheta_{t,u}:= x + \int_t^u \vartheta_r\dd S_r
$,
$t\le u\le T$, for $x\in\RR$, $t\in[0,T]$ and $\vartheta\in\Theta
$. (We
sometimes omit writing the dependence of $X^\vartheta$ on $t$.) This gives
%
\begin{equation}\label{equ2.22}
(X^\vartheta_u)^2
=x^2 + 2 \int_t^u X^\vartheta_{r-} \vartheta_r\dd S_r + \int_t^u
\vartheta_r^2 \dd[S]_r.
\end{equation}
Next we apply the product rule with (\ref{equ2.22}), (\ref{equ2.20}), (\ref{equ2.7}), (\ref{equ2.21}) and
then use $A = \int\lambda\dd\MM$ and
$
\q_-\SIdot[S] + [ \q, [S] ] = (\q_-+\Delta\q)\SIdot[S] = \q
\SIdot[S]
$
as well 
as (\ref{equ2.8}), (\ref{equ2.10}) for $\q$ and (\ref{equ2.15}) to obtain
%
\begin{eqnarray}
\label{equ2.23}\quad
(X^\vartheta_{t,u})^2 \q_u - x^2 \q_t
&=&
m_u - m_t
+ \int_t^u (X^\vartheta_{r-})^2 \dd\Bq_r\nonumber\\
&&{} + 2\int_t^u \q_{r-} X^\vartheta
_{r-} \vartheta_r \dd A_r
+\int_t^u \q_{r-} \vartheta_r^2 \dd[S]_r
\nonumber\\
&&{}
+ 2 \int_t^u X^\vartheta_{r-} \vartheta_r (\varphi_r + \lambda
_r\Delta\Bq_r ) \dd\MM_r + \int_t^u \vartheta_r^2\dd[ \q,[S] ]_r
\nonumber\\[-8pt]\\[-8pt]
&=&
m_u - m_t
+ \int_t^u (X^\vartheta_{r-})^2 \dd\Bq_r
\nonumber\\
&&{}+ \int_t^u \bigl( 2 X^\vartheta_{r-} \vartheta_r (
\varphi_r + \lambda_r \pq_r) +
\vartheta_r^2 \NNq{r} \bigr) \dd\MM_r
\nonumber\\
&=&
m_u - m_t
+ \int_t^u f(r,X^\vartheta_{t,r-}; \vartheta) \dd C_r,\nonumber
\end{eqnarray}
where $C\in\A^+_\loc(P)$ is a predictable process with
$\Bq= \int\beta\dd C$, $\MM= \int\nu\dd C$,
and
%
\begin{eqnarray}\label{equ2.24}
f(r,y;\vartheta)
&:=&
y^2 \beta_r + G_r(y, \vartheta_r) \nu_r\nonumber\\[-8pt]\\[-8pt]
&:=&
y^2 \beta_r + \bigl( 2y\vartheta_r (\varphi_r + \lambda_r \pq_r) +
\vartheta_r^2
\NNq{r} \bigr) \nu_r\nonumber
\end{eqnarray}
is a quadratic polynomial in $y$ with random processes as coefficients.
Replacing $C_t$ by $C_t+t$, we can assume that $C$, as well as its
continuous part~$C^c$, is strictly increasing.

By Corollary~\ref{Cor1.2}, $( (X^\vartheta_{t,u})^2 \q_u )_{t\le u\le
T}$ is a $P$-submartingale for every $\vartheta\in\Theta$ and a
$P$-martingale for the optimal $\thetast\in \Theta$, if that exists.
This means that the \mbox{$dC$-integral} in (\ref{equ2.23}) is
increasing for every $\vartheta\in\Theta$ and identically 0 for
$\vartheta= \thetast$, and the same then applies separately for the
corresponding integrals with respect to the continuous and purely
discontinuous parts $C^c$ and $C^d$ of $C$. Similarly as in
\citet{ManTev03N1}, we therefore obtain for each $x\in\RR$,
%
\begin{equation}\label{equ2.25}\qquad
\essinf_{\vartheta\in\Theta} f(r,x;\vartheta) = x^2 \beta_r + \nu_r
\essinf_{\vartheta\in\Theta} G_r(x, \vartheta _r) = 0,
\qquad\mbox{$P\otimes
C$-a.e.;}
\end{equation}
the details for this step are a bit more technical and are postponed to
step (2). Using the definition of $G_r(y,\vartheta_r)$ in (\ref{equ2.24}) and
completing the square gives
%
\begin{equation}\label{equ2.26}
G_r(x, \vartheta_r) = \NNq{r} \biggl( \vartheta_r + x {\varphi_r + \lambda_r
\pq_r \over \NNq{r} } \biggr)^2 - x^2 {(\varphi_r + \lambda_r \pq_r)^2
\over\NNq{r} },
\end{equation}
and we claim that for a localizing sequence $(\tau_n)_{n\in\NN}$,
%
\begin{equation}\label{equ2.27}
\vartheta^n:= -x {\varphi+ \lambda\pq\over{\mathcal N}(\q) } I_{\lsb0,
\tau _n\rsb} \in\Theta.
\end{equation}
Indeed, ${\mathcal N}(\q)$ is locally bounded away from 0 by Lemma~\ref{Lem2.3}, and
$ \pq$ is bounded like $\q$ due to Lemma~\ref{Lem1.5}. Moreover, $\int
\lambda^2 \dd\MM$ is locally bounded since it is predictable and RCLL,
and $\varphi$ is locally in $L^2(M)$ by construction. Thus we obtain via
Cauchy--Schwarz that both $\varphi$ and $\lambda$, and then also the ratio
in (\ref{equ2.27}), are locally in $L^2(M)\cap L^2(A) = \Theta$, as
claimed. Inserting $\vartheta^n$ into (\ref{equ2.26}) makes the first term in
(\ref{equ2.26}) vanish for $n\to\infty$ and thus yields
\[
\essinf_{\vartheta\in\Theta} G_r(x, \vartheta_r)
=
- x^2 {(\varphi_r + \lambda_r \pq_r)^2 \over\NNq{r} },
\qquad\mbox{$P\otimes C$-a.e.}
\]
Plugging this into (\ref{equ2.25}) and integrating gives
$
\Bq= \int\beta\dd C = \int{(\varphi+ \lambda\spP{\q})^2 \over
{\mathcal N}(\q)
} \dd\MM
$,
and plugging that in turn into (\ref{equ2.20}) shows that the triple $(\q,
\varphi, L^{\q})$ solves the BSDE (\ref{equ2.18}). Moreover, we see from
Lemma~\ref{Lem2.1} and $\q\le1$ that $\q$ is strictly positive and bounded.

(2) To prove (\ref{equ2.25}), we use the same basic approach as in \citet{ManTev03N1}, but we must be more careful and handle jumps since
$S$ is not
continuous. For ease of notation, we sometimes omit the third argument
$\vartheta$ of $f$. We first write $C = C^c + C^d$ and denote by
$(\tau
_k)_{k\in\NN}$ a sequence of stopping times exhausting the jumps of
$C^d$ (or $C$). Each $\tau_k$ is predictable because $C$ is
predictable. By Corollary~\ref{Cor1.2}, we then have with probability 1
that $C_{\fatdot}(\omega)$ is RCLL and simultaneously for all rational
$s\in[0,T]$ that
%
\begin{eqnarray}
\label{equ2.28}
&&\int_s^u f(r, X^\vartheta_{s, r-}; \vartheta) \dd C_r,
\qquad s\le u\le T,\qquad \mbox{is increasing},
\\
&&
\label{equ2.29}
\int_s^u f(r, X^\vartheta_{s, r-}; \vartheta) \dd C^c_r,
\qquad s\le u\le T,\qquad \mbox{is increasing,}
\end{eqnarray}
for each $\vartheta\in\Theta$, and for the optimal $\vartheta
^{*,s}$, the
processes in (\ref{equ2.28}) and (\ref{equ2.29}) vanish identically. Indeed,
(\ref{equ2.29}) follows from (\ref{equ2.28}) since the process in (\ref{equ2.29}) is simply the
continuous part of the process in (\ref{equ2.28}). For
any $\tau_k$, we thus have with probability 1 that
\[
\int_s^{\tau_k(\omega)} f(r,X^\vartheta_{s, r-}; \vartheta)
(\omega) \dd C_r(\omega) \ge0
\qquad\mbox{for all rational }s<\tau_k(\omega).
\]
Because $\tau_k$ is predictable, there are stopping times $(\sigma
^{(n)}_k)_{n\in\NN}$ taking only rational values and such that $\lim
_{n\to\infty} \sigma^{(n)}_k = \tau_k$ and $\sigma^{(n)}_k < \tau
_k$ on
$\{\tau_k>0\} = \Omega$; see Theorem IV.77 in \citet{DelMey78}. Thus we obtain for
$P$-almost all $\omega$ that
\[
\int_{\sigma^{(n)}_k(\omega)}^{\tau_k(\omega)} f\bigl( r,
X^\vartheta
_{\sigma^{(n)}_k, r-}; \vartheta\bigr) (\omega) \dd C_r(\omega) \ge0
\qquad\mbox{for all $k$ and $n$.}
\]
These integrals tend to
$
f(\tau_k, X^\vartheta_{\tau_k-, \tau_k-}; \vartheta) (\omega)
\Delta C_{\tau
_k} (\omega)
=
f(\tau_k,x;\vartheta) (\omega) \Delta C_{\tau_k}
(\omega)
$ as $n\to\infty$
because $X^\vartheta_{\tau_k-, \tau_k-} = x$, and so we get
%
\begin{equation}\label{equ2.30}
f(\tau_k,x;\vartheta) \Delta C_{\tau_k} \ge0
\qquad\mbox{for all $k\in\NN$,
$P$-a.s.,}
\end{equation}
which means that $f( \cdot, x; \vartheta) \ge0$ $P\otimes C^d$-a.e.,
for each $\vartheta\in\Theta$. For the optimal $\vartheta^{*,s}$,
we get the
null process in (\ref{equ2.28}), hence equality in (\ref{equ2.30}), and so we have
%
\begin{equation}\label{equ2.31}
\essinf_{\vartheta\in\Theta} f( \cdot,x;\vartheta) = 0,
\qquad\mbox{$P\otimes
C^d$-a.e.}
\end{equation}

For the continuous part $C^c$, (\ref{equ2.29}) gives with
$
\tau_s(\varepsilon):= \inf\{ t\ge s | C^c_t \ge C^c_s +
\varepsilon\}
$
that
%
\begin{equation}\label{equ2.32}\quad
\int_s^{\tau_s(\varepsilon)} f(t, X^\vartheta_{s, t-};
\vartheta) \dd C^c_t \ge0\qquad\mbox{for all
rational $s\in[0,T]$, $P$-a.s.}
\end{equation}
We claim that for each $u\ge s$,
%
\begin{equation}\label{equ2.33}
s\mapsto\int_s^u f(t, X^\vartheta_{s, t-}; \vartheta) \dd C^c_t
\qquad\mbox{is $P$-a.s. right-continuous.}
\end{equation}
Postponing the argument for the moment, we obtain that the inequality
in (\ref{equ2.32}) also holds for all $s\in[0,T]$, $P$-a.s. Setting
$
\sigma_t(\varepsilon):= \inf\{ s\ge0 | C^c_s \ge C^c_t -
\varepsilon\}
$,
we then get as in Appendix B of \citet{ManTev03N1} via
Fubini's theorem that
(dropping arguments $\vartheta$ from $f$)
%
\begin{eqnarray}
\label{equ2.34}
&&
\int_0^T \biggl| {1\over\varepsilon} \int_s^{\tau_s(\varepsilon)}
f(t, X^\vartheta
_{s, t-}) \dd C^c_t - f(s,x) \biggr| \dd C^c_s\nonumber\\
&&\qquad \le
\int_0^T {1\over\varepsilon} \int_{\sigma_t(\varepsilon)}^t |
f(t, X^\vartheta
_{s, t-}) - f(t,x) | \dd C^c_s \dd C^c_t
\\
&&\qquad\quad {}
+ \int_0^T {1\over\varepsilon} \int_s^{\tau
_s(\varepsilon)} | f(t, x) -
f(s,x) | \dd C^c_t \dd C^c_s;
\nonumber
\end{eqnarray}
this uses that $C^c_{\tau_s(\varepsilon)} - C^c_s = \varepsilon$ by
continuity of $C^c$. The second term on the right-hand side of (\ref{equ2.34})
tends to 0 as $\varepsilon\searrow0$ by Corollary B.1 in \citet{ManTev03N1}. Writing
\begin{eqnarray*}
b^\varepsilon_t
:\!&=&
\sup\{ |X^\vartheta_{s, t-} - x| | \sigma
_t(\varepsilon) < s < t
\}\\
&=&
\sup\biggl\{ \biggl| \int_s^{t-} \vartheta_r \dd S_r \biggr|\bigg|
\sigma
_t(\varepsilon) < s < t \biggr\},
\end{eqnarray*}
we have $\sigma_t(\varepsilon)\nearrow t$ for $\varepsilon\searrow
0$ by
continuity of $C^c$ and therefore $b^\varepsilon_t\searrow0$ as
$\varepsilon
\searrow0$. Moreover, we have (uniformly in $\varepsilon$ and $t$)
$
b^\varepsilon_t \le2 \sup_{0\le r\le T} |\vartheta\SIdot S_r|
$
which is in $L^2(P)$, hence $P$-a.s. finite, for $\vartheta\in\Theta$. The
first term on the right-hand side of (\ref{equ2.34}) can now be estimated
above by
\[
\int_0^T \sup\bigl\{ |f(t,y; \vartheta) - f(t,x; \vartheta)|
\big|
|y-x|\le b^\varepsilon_t \bigr\} \dd C^c_t =: \int_0^T h_\varepsilon
(t; \vartheta
) \dd C^c_t
\]
since $C^c_t - C^c_{\sigma_t(\varepsilon)} = \varepsilon$ by
continuity of
$C^c$. Now we use the definition of $f$ in (\ref{equ2.24}) to obtain
\[
h_\varepsilon(t; \vartheta)
\le
(b^\varepsilon_t)^2 |\beta_t| + b^\varepsilon_t \bigl( 2 |\beta_t|
|x| + 2 \nu
_t |\vartheta_t| (|\varphi_t| + |\lambda_t| \pq_t) \bigr).
\]
This shows that $P$-a.s., $h_\varepsilon(t; \vartheta)\to0$ for all $t$ as
$\varepsilon\searrow0$. Moreover, $b^\varepsilon_t$ can be bounded uniformly
in $\varepsilon$ and $t$, $P$-a.s., and using
%
\begin{eqnarray}
\label{equ2.35}
&&\int_0^T |\beta_t| \dd C^c_t
\le
\int_0^T |d\Bq_t|,
\\
\label{equ2.36}
&&\int_0^T \nu_t |\vartheta_t| (|\varphi_t| +
|\lambda_t| \pq_t) \dd C^c_t\nonumber\\[-8pt]\\[-8pt]
&&\qquad\le
\biggl( \int_0^T \vartheta_t^2 \dd\MM_t \biggr)^{1/2} \biggl( 2
\int_0^T
\biggl(\varphi_t^2 + \lambda_t^2\biggr) \dd\MM_t \biggr)^{1/2}\nonumber
\end{eqnarray}
shows that we can apply dominated convergence to get
$
\int_0^T h_\varepsilon(t; \vartheta) \dd C^c_t \longrightarrow0
$
as $\varepsilon\searrow0$, $P$-a.s. With a similar argument, we can prove
(\ref{equ2.33}). Indeed, for $s_n\searrow s$, we have
\begin{eqnarray*}
&&
\biggl| \int_s^u f(t, X^\vartheta_{s,t-}) \dd C^c_t - \int_{s_n}^u f(t,
X^\vartheta_{s_n,t-}) \dd C^c_t \biggr|\\
&&\qquad \le
\int_s^{s_n} |f(t, X^\vartheta_{s,t-})| \dd C^c_t
+ \int_{s_n}^u |f(t, X^\vartheta_{s,t-}) - f(t, X^\vartheta
_{s_n,t-})| \dd C^c_t,
\end{eqnarray*}
and the first term on the right-hand side tends to 0 $P$-a.s. as $n\to
\infty
$ by continuity of $C^c$. Writing
$
h_n(t):= |f(t, X^\vartheta_{s,t-}) - f(t, X^\vartheta_{s_n,t-})|
$,
we have $h_n(t)\to0$ as $n\to\infty$ by the right-continuity of the
stochastic integral and since $f$ from (\ref{equ2.24}) is continuous with
respect to the second argument $y$. So (\ref{equ2.33}) will follow by
dominated convergence as soon as we show that
%
\begin{equation}\label{equ2.37}
\int_0^T \sup_{n\in\NN} h_n(t) \dd C^c_t < \infty,
\qquad\mbox{$P$-a.s.}
\end{equation}
But the definition of $f$ in (\ref{equ2.24}) yields that
\[
h_n(t)
\le
4 |\beta_t| \Bigl( |x|^2 + \sup_{0\le r\le T} |\vartheta\SIdot
S_r|^2 \Bigr)
+ 2 |\nu_t| |\vartheta_t| ( |\varphi_t| + |\lambda_t| \pq_t) \sup
_{0\le r\le
T} |\vartheta\SIdot S_r|,
\]
and so (\ref{equ2.37}) follows again by (\ref{equ2.35}) and (\ref{equ2.36})
because $\sup_{0\le r\le T} |\vartheta\SIdot S_r| < \infty$ $P$-a.s. This
establishes (\ref{equ2.33}).

Putting together all the results so far, (\ref{equ2.34}) therefore yields
that with probability 1, we have
$
{1\over\varepsilon} \int_s^{\tau_s(\varepsilon; \vartheta)} f(t,
X^\vartheta_{s,
t-}; \vartheta) \dd C^c_t \longrightarrow f(s,x; \vartheta)
$
in $L^1(dC^c)$ as $\varepsilon\searrow0$. Together with (\ref{equ2.32}), this
gives $f( \cdot, x;\vartheta) \ge0$ $P\otimes C^c$-a.e., for each
$\vartheta\in\Theta$. For the optimal $\vartheta^{*,s}$, we again get
equality so that finally
\[
\essinf_{\vartheta\in\Theta} f( \cdot, x;\vartheta) = 0,
\qquad\mbox{$P\otimes C^c$-a.e.,}
\]
and combining this with (\ref{equ2.31}) yields (\ref{equ2.25}).

(3) We next show that $\thetast$ for fixed $t$ is given by (\ref{equ2.19}).
Since $(\q, \varphi, L^{\q})$ satisfies (\ref{equ2.18}), It\^o's formula
gives via (\ref{equ2.22}) and (\ref{equ2.8})--(\ref{equ2.11}) like in (\ref{equ2.23}) for any $\vartheta
\in\Theta$ that
%
\begin{eqnarray}
\label{equ2.38}
&&(X^\vartheta_u)^2 \q_u - x^2 \q_t\nonumber\\
&&\qquad=
m_u - m_t\nonumber\\
&&\qquad\quad{}+ \int_t^u \biggl( (X^\vartheta_{r-})^2 {(\varphi_r
+ \lambda_r \pq_r)^2 \over
\NNq{r} } + 2 \q_{r-} X^\vartheta_{r-} \vartheta_r \lambda_r
\nonumber\\[-8pt]\\[-8pt]
&&\hspace*{33pt}\qquad\quad{}
+ \q_{r-} \vartheta_r^2 (1+\lambda_r^2 \Delta\MM_r)
+ 2 X^\vartheta_{r-} \vartheta_r (\varphi_r + \lambda_r \Delta\Bq_r)
\nonumber\\
&&\hspace*{95.4pt}\qquad\quad{}
+ \vartheta_r^2 \bigl( \Delta\Bq_r (1+\lambda_r^2 \Delta\MM_r) +
\func_r(\q
) \bigr) \biggr)\dd\MM_r
\nonumber\\
&&\qquad=
m_u - m_t
+ \int_t^u \biggl( \vartheta_r \sqrt{\NNq{r}} + X^\vartheta_{r-}
{\varphi_r +
\lambda_r \pq_r \over\sqrt{\NNq{r}} }\biggr)^2 \dd\MM_r.
\nonumber
\end{eqnarray}
By Corollary~\ref{Cor1.2}, the process in (\ref{equ2.38}) is a martingale on
$\lsb t,T\rsb$ for the optimal $\thetast$, and so
%
\begin{equation}\label{equ2.39}
\thetast= - X^{\thetast}_- {\varphi+ \lambda\pq\over{\mathcal N}(\q) },
\qquad\mbox{$P\otimes\MM$-a.e. on $\rsb t,T\rsb$.}
\end{equation}
Integrating with respect to $S$ thus shows for $x=1$ that $X^{\thetast}
= 1 + \int_t^{\fatdot} \thetast\dd S$ satisfies
the linear SDE
$
X^{\thetast}_u = 1 - \int_t^u X^{\thetast}_{r-} {\varphi_r +
\lambda_r \spP
{\q}_r \over\NNq{r} } \dd S_r
$
for $t\le u\le T$, and this implies that
%
\begin{equation}\label{equ2.40}
X^{\thetast} = \tE\biggl( -{\varphi+ \lambda\pq\over{\mathcal N}(\q) }
\SIdot S\biggr).
\end{equation}
Because $\thetast$ is in $\Theta$, we have $X^{\thetast}\in\SS
^2(P)$ so
that the stochastic exponential is indeed in
$\SS^2(P)$; and plugging (\ref{equ2.40}) into (\ref{equ2.39}) yields the
expression (\ref{equ2.19}) for~$\thetast$. Since $t$ was arbitrary, we
have now shown that (a) implies (b) and that we then have (\ref{equ2.19}).

(4) Conversely, let us start from (b). Again fix $t$. Using the fact that
$(Y,\psi,L)$ solves the BSDE (\ref{equ2.18}),
we obtain completely analogously as for (\ref{equ2.38}) for any $\vartheta\in
L(S)$ that
%
\begin{eqnarray}\label{equ2.41}
&&(X^\vartheta_u)^2 Y_u - x^2 Y_t\nonumber\\[-8pt]\\[-8pt]
&&\qquad =
m_u - m_t
+ \int_t^u \biggl( \vartheta_r \sqrt{\NNY{r}} + X^\vartheta_{r-}
{\psi_r +
\lambda_r \pY_r \over\sqrt{\NNY{r}} } \biggr)^2 \dd\MM_r
\nonumber
\end{eqnarray}
for $t\le u\le T$. So $(X^\vartheta)^2 Y$ is a local $P$-submartingale on
$\lsb t,T\rsb$; but since $Y$ is bounded and
$1+\vartheta\SIdot S\in\SS^2(P)$ for $\vartheta\in\Theta$, we get that
$(X^\vartheta)^2 Y$ is actually a true $P$-submartingale on $\lsb
t,T\rsb$
so that $Y_T=1$ gives
$
Y_t
\le
E[ ( 1 + \int_t^T \vartheta_r\dd S_r )^2 | \F
_t ]
$
for any $\vartheta\in\Theta$. The definition\vspace*{1pt} in (\ref{equ2.1}) thus yields
$
Y_t \le\Vn_t(1) = \q_t
$
for all $t\in[0,T]$. To prove the converse inequality, define the
predictable process $\thetat^{(t)}$ by the right-hand side of (\ref{equ2.19}).
Integrating then shows as for (\ref{equ2.40}) that
\mbox{$
X^{\thetat^{(t)}} = \tE( - {\psi+ \lambda\spP{Y} \over{\mathcal N}(Y) }
\SIdot S),
$}
and because this stochastic exponential is in $\SS^2(P)$ by the
assumption in b), we see that $\thetat^{(t)}$ coming from (\ref{equ2.19})
is actually in $\Theta$. Plugging $\thetat^{(t)}$ into (\ref{equ2.41})
shows by (\ref{equ2.19}) that the $d\MM$-integral vanishes; so $(X^{\thetat
^{(t)}})^2 Y$ is a $P$-martingale on $\lsb t,T\rsb$ and hence
$
Y_t
=
E[ ( 1 + \int_t^T \thetat^{(t)}_r\dd S_r )^2
| \F_t
]
\ge
\Vn_t(1)
=
\q_t
$
by (\ref{equ2.1}). So we obtain \mbox{$Y=\q$}, hence also $\psi\SIdot M = \varphi
\SIdot M$, $L=L^{\q}$, and this shows that any
solution of (\ref{equ2.18}) with the properties in (b) coincides with $(\q,
\varphi, L^{\q})$, giving uniqueness. Finally, $Y=\q$
shows that $(X^\vartheta)^2 \q$ is a\vspace*{1pt} $P$-submartingale on $\lsb
t,T\rsb$
for any $\vartheta\in\Theta$ and a $P$-martingale for
$\vartheta=\thetat^{(t)}\in\Theta$; so $\thetat^{(t)}$ is optimal by
Corollary~\ref{Cor1.2} and in particular, an optimal
$\thetast(1,0) = \thetat^{(t)}$ exists. Since $t$ was arbitrary, we
have also shown that (b) implies (a), and part (1) of
Theorem~\ref{Theorem2.4} is proved.

(5) It remains to prove part (2). But if there is some $Q\in\PPes$ with
$R_2(P)$, the space $ L^2(\F_t,P)+G_{t,T}(\Theta) =\{ X +
\vartheta\SIdot S_T | X\in L^2(\F _t,P), \vartheta\in\Theta_{t,T}\} $
is closed in $L^2(P)$ by Theorem 5.2 of \citet{ChoKraStr98}, for
every $t$, so that an optimal $\thetast$ exists. Moreover, we then have
$\q\ge\delta>0$ by Lemma~\ref{Lem2.1}, and so part (2) follows directly
from part (1).
\end{pf}
\begin{Rem*}
If $d>1$, the backward equation (\ref{equ2.18}) looks more
complicated. Using the notation from the remark before Theorem
\ref{Theorem2.4}, in particular (\ref{equ2.16}) and (\ref{equ2.17}), the equation reads
\begin{eqnarray*}
Y_t
&=&
Y_0 + \int_0^t {( \psi_s + \lambda_s \pY_s)\tr\mu_s ( \NNY
{s}
)^{-1} \mu_s ( \psi_s + \lambda_s \pY_s) } \dd B_s\\
&&{} +
\int_0^t \psi_s\dd M_s + L_t,
\\
Y_T&=&1,
\end{eqnarray*}
where
$
\NNY{s}:= \pY_s ( \mu_s + (\mu_s\lambda_s)\tr\mu_s\lambda_s
\Delta
B_s ) + \func_s(Y)
$.
We do not give details.
\end{Rem*}

For later use, we record the following consequence of Theorem~\ref{Theorem2.4}.
%
\begin{corollary}\label{Cor2.5}
Under the assumptions of Theorem~\ref{Theorem2.4}, suppose \textup{(a)}
or \textup{(b)} there hold. Define
%
\begin{equation}\label{equ2.42}
\gamma:= - {\psi+ \lambda\pY\over{\mathcal N}(Y) }
=
- {\psi+ \lambda\pY\over\pY(1+\lambda^2\Delta\MM) + \func
(Y) },
\end{equation}
where $(Y, \psi, L)$ is the solution of the BSDE (\ref{equ2.18}), and
recall the process $\vone$ from the quadratic representation (\ref{equ1.6}) of
$\VH$. For every $t\in[0,T]$, we then have
%
\begin{equation}\label{equ2.43}
\vone_t = E[ H \tE(\gamma\SIdot S)_T | \F_t ],
\qquad\mbox{$P$-a.s.}
\end{equation}
and the process $( \tE(\gamma\SIdot S)_u \vone_u )_{t\le u\le T}$ is
a $P$-martingale on $\lsb t,T\rsb$.
\end{corollary}
\begin{pf}
Fix $t$. Because we have
$
1 + \int_t^T \thetast_r(1,0) \dd S_r
=
X^{\thetast}_T
=
\tE(\gamma\SIdot S)_T
$
by (\ref{equ2.40}) and the definition (\ref{equ2.42}) of $\gamma$, (\ref{equ2.43}) follows
directly from (\ref{equ1.8}). Moreover, it is easy to check
that for any semimartingale $X$ and any $u\le T$, we have
$
\uE(X)_T = {\E(X)_T \over\E(X)_u}
$
$P$-a.s. on $\{\E(X)_u\ne0\}$
and $\E(X)_T=0$ $P$-a.s. on $\{\E(X)_ u =0\}$. Taking $X:= \gamma
\SIdot S
- (\gamma\SIdot S)^t$, $u\ge t$ and setting for brevity
$
D_ u:=\{ \tE(\gamma\SIdot S)_ u \ne0\}
$
therefore gives the desired martingale property via
\begin{eqnarray*}
\tE(\gamma\SIdot S)_ u \vone_ u
&=&
I_{D_ u} \tE(\gamma\SIdot S)_ u E[ H \uE(\gamma\SIdot
S)_T
| \F_ u ]
\\
&=&
I_{D_ u} E[ H \tE(\gamma\SIdot S)_T | \F_ u ]
\\
&=&
E[ H \tE(\gamma\SIdot S)_T | \F_ u ];
\end{eqnarray*}
integrability holds since $H\in L^2(P)$ and $\tE(\gamma\SIdot S)\in
\SS
^2(P)$ by part (1b) of Theorem~\ref{Theorem2.4}.
\end{pf}

As before, we can connect our results to the dual problem, as follows.
%
\begin{proposition}\label{Prop2.6}
Under the assumptions\vspace*{1pt} of Theorem~\ref{Theorem2.4}, suppose
\textup{(a)} or \textup{(b)} there hold.
Then the variance-optimal signed martingale measure\vadjust{\goodbreak} $\Qt\in\PPss$ is
given by
%
\begin{equation}\label{equ2.44}
{d\Qt\over dP} = {1\over Y_0} \E\biggl( - {\psi+\lambda\pY\over{\mathcal
N}(Y) } \SIdot S \biggr)_T = {1\over Y_0} \E( \gamma\SIdot S)_T,
\end{equation}
where $(Y,\psi,L)$ is the solution of the BSDE (\ref{equ2.18}). If we have,
in addition, that
%
\begin{eqnarray}\label{equ2.45}
\gamma_t \Delta S_t
&=&
- {\psi_t + \lambda_t \pY_t \over\pY_t (1+\lambda_t^2 \Delta\MM
_t) +
\func_t(Y) } \Delta S_t\nonumber\\[-8pt]\\[-8pt]
&>& -1,
\qquad\mbox{$P$-a.s. for $0\le t\le T$,}\nonumber
\end{eqnarray}
then the VOMM exists and is given by $\Qt$ from (\ref{equ2.44}).
\end{proposition}
\begin{pf}
From the BSDE (\ref{equ2.18}) and It\^o's formula, we obtain by
straightforward computation that the product $Y
\E( - {\psi+\lambda\spP{Y} \over{\mathcal N}(Y)} \SIdot S)$ is a local
$P$-martingale, and it is even a true $P$-martingale since $Y$ is
bounded and the stochastic exponential is in $\SS^2(P)$, and so (\ref{equ2.44})
defines a signed measure $\Qt\ll P$ with $P$-square-integrable
density process $Z^{\Qt} = {Y \over Y_0} \E( - {\psi+\lambda\spP{Y}
\over{\mathcal N}(Y) } \SIdot S)$ and $\Qt[\Omega]=1$. Note
for\vspace*{1pt}
(\ref{equ2.44}) that $Y_T=1$. Another straightforward but slightly lengthier
computation shows that $Z^{\Qt} S$ is a local $P$-martingale so that
$\Qt\in\PPss$. Finally, the representation (\ref{equ2.44}) of $d\Qt\over
dP$ as a constant plus a ``good'' stochastic integral of $S$
implies that $\Qt$ is variance-optimal; see, for instance, Lemma~2.1 in
\citet{DelSch96}. Note here that the same argument as
in step\vspace*{1pt} (4) of the proof
of Theorem~\ref{Theorem2.4} implies that the integrand $\vartheta:= {1\over Y_0}
\gamma\E(\gamma\SIdot S)_-$ is\vspace*{-1pt} in $\Theta$ so that $\vartheta
\SIdot S$ is
a $Q$-martingale for every $Q\in\PPes$. If (\ref{equ2.45}) holds, then
clearly $Z^{\Qt}>0$; so $\Qt$ is then equivalent to $P$, hence in
$\PPes$,
and is the VOMM.
\end{pf}
\begin{Rem*}
From (\ref{equ2.43}), the proof of Proposition~\ref{Prop2.6} and
$Y=\vtwo$, we can see that under the assumptions of Theorem~\ref{Theorem2.4}
and (\ref{equ2.45}), the process
\[
\vone\E(\gamma\SIdot S) = \vone Y_0
Z^{\Qt} / Y
\]
is a $P$-martingale with final value $H \E(\gamma
\SIdot S)_T = H Y_0 Z^{\Qt}_T$. This implies that
\[
{\vone_t \over\vtwo_t}
=
{\vone_t \over Y_t}
=
E_{\Qt} [ H | \F_t ],
\qquad 0\le t\le T.
\]
\end{Rem*}

%
%
\section{Mean-variance hedging: From $(1,0)$ to $(x,H)$}\label{sec3}

Recall from Theorem~\ref{Theorem1.4} that the dynamic value process of the
mean-variance hedging problem has the quadratic form
\[
\VH(x) = \vnull- 2\vone x + \vtwo x^2.
\]
Our goals in this section are to describe the coefficient processes
$\vnull, \vone,\break \vtwo$ via backward stochastic differential equations
(BSDEs) and to give explicit expressions for the optimal strategies
$\thetast(x,H)$. This will be done under the same assumptions as in Section
\ref{sec2}.

A general solution for the MVH problem has been given by
\citet{CerKal07} in their
Theorem 4.10 and Corollary 4.11. However, that solution
involves either a process $N$ which is very hard to find [see
\citet{CerKal07},
Definition 3.12] or the variance-optimal martingale measure
[called $Q^*$ in \citet{CerKal07}; see their Proposition
3.13] which is also
notoriously difficult to determine. With our approach, we can be more explicit.

To formulate our main result, we introduce the system of BSDEs,
%
\begin{eqnarray}
\label{equ3.1}\qquad
d\Ytwo_s
&=&
{(\psitwo_s + \lambda_s {}^{\p}\Ytwo_s)^2 \over\N
_s(\Ytwo)} \dd\MM_s +
\psitwo_s\dd M_s + d\Ltwo_s,\qquad
\Ytwo_T=1,
\\
\label{equ3.2}
d\Yone_s
&=&
{(\psitwo_s+\lambda_s {}^{\p}\Ytwo_s) (\psione
_s+\lambda_s {}^{\p}\Yone_s)
\over\N_s(\Ytwo)} \dd\MM_s\nonumber\\[-8pt]\\[-8pt]
&&{}+ \psione_s\dd M_s + d\Lone_s,\qquad
\Yone_T=H,\nonumber
\\
\label{equ3.3}
d\Ynull_s
&=&
{(\psione_s+\lambda_s {}^{\p}\Yone_s)^2 \over\N
_s(\Ytwo)} \dd\MM_s +
d\Nnull_s,\qquad
\Ynull_T = H^2.
\end{eqnarray}
A \textit{solution} of this system consists of tuples $(\Ytwo, \psitwo,
\Ltwo)$, $(\Yone, \psione, \Lone)$, $(\Ynull, \Nnull)$ where $\psitwo,
\psione$ are in $L^1_\loc(M)$; $\Ltwo, \Lone$ are in $\M_{0,\loc}(P)$
and strongly $P$-orthogonal to $M$; $\Nnull$ is a local $P$-martingale;
and $\Ytwo, \Yone, \Ynull$ are $P$-special semimartingales with $[
N^{\Ytwo}, [S] ] \in\A_\loc(P)$. We point out that (\ref{equ3.1}) is
the same equation as (\ref{equ2.18}) before Theorem~\ref{Theorem2.4}.
Note also that (given $\Ytwo, \psitwo, \Ltwo$) the equation
(\ref{equ3.2}) is linear and can therefore be solved explicitly; and
$\Ynull$ and $\Nnull $ for (\ref{equ3.3}) can even be written down
directly. In the case where $S$ is continuous, this system has been
obtained and studied in \citet{ManTev03N1} or (under the
additional assumption that $\FF$ is continuous) in
\citet{BobSch04}. For a Markovian setting within a Brownian
filtration, the corresponding PDEs can also be found in
\citet{BerKogLo01}, with a heuristic treatment.
%
\begin{theorem}\label{Theorem3.1}
Suppose (as in Theorem~\ref{Theorem2.4}) that
$S\in\SS^2_\loc(P)$ and\break $\PPes\ne\varnothing$, and fix $H\in L^2(\F_T,
P)$. Then:

\begin{enumerate}[(2)]
\item[(1)] The following two assertions are equivalent:

\begin{enumerate}[(a)]
\item[(a)] For every $t\in[0,T]$, there exists an optimal
$\thetast(x,H)\in\Theta_{t,T}(0)$ for (\ref{equ1.2}) for every $x\in\RR$.

\item[(b)] For each $x\in\RR$, there is a solution to the BSDE
system (\ref{equ3.1})--(\ref{equ3.3}) with:

\begin{enumerate}[(iii)]
\item[(i)] $\Ltwo\in\M^2_{0,\loc}(P)$, $\psitwo\in L^2_\loc(M)$,
$\Ytwo$ bounded and strictly positive, and with the
property that for every $t\in [0,T]$, the process $( \tE( -
{\psitwo+ \lambda\spP{\Ytwo} \over\N(\Ytwo)} \SIdot S)_u )_{t\le u\le
T}$ is in $\SS^2(P)$;\vspace*{2pt}

\item[(ii)] $\Lone\in\M^2_{0,\loc}(P)$, $\psione\in L^2_\loc(M)$,
$|\Yone|^2$ of class (D), and such that for every $t\in[0,T]$, the
solution $X^{(t)}$ of the linear SDE
%
\begin{equation}\label{equ3.4}
X^{(t)}_u
=
x + \int_t^u {\psione_r + \lambda_r {}^{\p}\Yone
_r \over\N_r(\Ytwo)} \dd S_r
- \int_t^u {\psitwo_r + \lambda_r {}^{\p}\Ytwo_r
\over\N_r(\Ytwo)}
X^{(t)}_{r-} \dd S_r
\end{equation}
on $\lsb t,T\rsb$ is in $\SS^2(P)$;

\item[(iii)] $\Ynull$ is a true $P$-submartingale and (hence) of
class (D).
\end{enumerate}
\end{enumerate}
If \textup{(a)} or \textup{(b)} holds, then the value process $\VH$
from (\ref{equ1.2}) admits the representation
%
\begin{equation}\label{equ3.5}
\VH(x) = \vnull- 2\vone x + \vtwo x^2,
\end{equation}
where the processes $\vtwo, \vone, \vnull$ satisfy the BSDE system
(\ref{equ3.1})--(\ref{equ3.3}), and for every $t\in[0,T]$, the optimal
wealth process
$
X^{\thetast}_u = x + \int_t^u \thetast_r(x,\break H) \dd S_r
$,
$t\le u\le T$, satisfies the SDE (\ref{equ3.4}) and $\thetast= \thetast
(x,H)$ is given by the feedback formula
%
\begin{equation}\label{equ3.6}
\thetast_u
=
{\psione_u + \lambda_u {}^{\p}\Yone_u \over\N
_u(\Ytwo)}
- {\psitwo_u + \lambda_u {}^{\p}\Ytwo_u \over\N
_u(\Ytwo)} X^{\thetast}_{u-},
\qquad t\le u\le T.
\end{equation}

\item[(2)] Suppose, in addition, that there is some $Q\in
\PPes$ satisfying
the reverse H\"older inequality $R_2(P)$. Then the value process $\VH$
from (\ref{equ1.2}) has the form (\ref{equ3.5}), where the processes $\vtwo
, \vone, \vnull$ are those unique solutions of the BSDE system
(\ref{equ3.1})--(\ref{equ3.3}) for which $\Ynull$ and $|\Yone|^2$ are of
class (D) and $c\le\Ytwo\le C$ for constants $C\ge c>0$. Moreover, for
every $t\in[0,T]$, the optimal strategy $\thetast(x,H)$ for (\ref{equ1.2})
exists, and its wealth process $X^{\thetast}$ satisfies the SDE
(\ref{equ3.4}).
\end{enumerate}
\end{theorem}
\begin{Rem*}
The integrability condition on the exponential in (i) is not really
needed. In fact, like in the proof of Theorem~\ref{Theorem1.4}, one can
argue that we have $\thetast(1,0) = \thetast(1,H) - \thetast(0,H)$ so that the
integrability required in (i) follows from that in~(ii). But for
simpler comparison with Theorem~\ref{Theorem2.4}, we have kept the
formulation as a condition.
\end{Rem*}
\begin{pf*}{Proof of Theorem~\ref{Theorem3.1}}
As in the proof of Theorem~\ref{Theorem2.4}, we denote by $m$ a generic local
$P$-martingale that can change from one appearance to the
next.\vadjust{\goodbreak}

(1) We first note that as in Theorem~\ref{Theorem1.4}, the existence of
optimal strategies $\thetast(1,0)$ (for $x=1, H\equiv0$) follows from
(a) and is, by Theorem~\ref{Theorem2.4}, equivalent to the solvability
of (\ref{equ3.1}) such that (i) holds in (b). So let us start from (a).
We note that (\ref{equ3.5}) holds due to Theorem~\ref{Theorem1.4}, and
first derive the BSDE for $\vone$. By Lem\-ma~\ref{Lem1.5} and the
Galtchouk--Kunita--Watanabe decomposition, we have
%
\begin{equation}\label{equ3.7}
\vone= \vone_0 + \mone+ \aone= \vone_0 + \psione\SIdot M + \Lone+
\aone
\end{equation}
with $\psione\in L^2_\loc(M)$, $\Lone\in\M^2_{0,\loc}(P)$ strongly
$P$-orthogonal to $M$, and $\aone$ predictable and of finite variation.
Exactly as for (\ref{equ2.21}), this yields
%
\begin{equation}\label{equ3.8}
\bigl[\vone, S\bigr] = m + \bigl(\psione+ \lambda\Delta\aone\bigr)\SIdot\MM.
\end{equation}
Now fix $t$, recall $\gamma$ from (\ref{equ2.42}) in Corollary~\ref{Cor2.5}
and write $\E:= \tE(\gamma\SIdot S)$ for brevity. Then combining
$d\E
= \E_- \gamma\dd S$ with the product rule, (\ref{equ3.7}), (\ref{equ2.7}),
(\ref{equ3.8}) and (\ref{equ2.10}) yield
%
\begin{eqnarray}
\label{equ3.9}
\E\vone
&=&
m + \E_- \SIdot\aone+ \bigl(\vone_- \E_- \gamma\lambda\bigr)\SIdot\MM\nonumber\\
&&{}+ \bigl( \E_- \gamma\bigl(\psione+ \lambda\Delta\aone\bigr)
\bigr)\SIdot\MM
\\
&=&
m + \E_- \SIdot\bigl( \aone+ \bigl( \gamma\bigl(\psione+ \lambda
{}^{\p}\vone\bigr) \bigr) \SIdot\MM\bigr).
\nonumber
\end{eqnarray}
But we know from Corollary~\ref{Cor2.5} that $\E\vone$ is a $P$-martingale
on $\lsb t,T\rsb$, and so the predictable finite variation term on the
right-hand side of (\ref{equ3.9}) must be identically zero. With
$C\in\A^+_\loc(P)$ predictable and such that $\aone\ll C$, $\MM
\ll
C$, we thus obtain that the process
$
\int\tE(\gamma\SIdot S)_- \{ {d\aone\over dC} + \gamma(\psione+
\lambda{}^{\p}\vone) \,{d\MM\over dC} \} \dd C
$
vanishes identically. Since $\tE(\gamma\SIdot S)_t = 1$, we can argue
analogously to steps (1) and (2) in the proof of Theorem~\ref{Theorem2.4} to get
\[
{d\aone\over dC} + \gamma\bigl(\psione+ \lambda{}^{\p}\vone\bigr)\, {d\MM\over
dC} = 0,
\qquad\mbox{$P\otimes C$-a.e.}
\]
Integrating with respect to $C$ gives
\[
\aone
=
{
- \int\gamma\bigl(\psione+ \lambda{}^{\p}\vone\bigr)\dd\MM
}
=
\int{(\psitwo+ \lambda{}^{\p}\Ytwo) (\psione+
\lambda{}^{\p}\vone) \over\N(\Ytwo)} \dd\MM,
\]
and plugging this into (\ref{equ3.7}) shows that $(\vone, \psione, \Lone
)$ satisfies the BSDE (\ref{equ3.2}). Moreover, as already used, we know
from Lemma~\ref{Lem1.5} that $|\vone|^2$ is of class (D), and it only
remains for (ii) to check the last integrability property.

(2) We next argue that the BSDE (\ref{equ3.3}) has a solution, starting
with a calculation that is used again later. Fix $t$, take any
$\vartheta$
in $\Theta$ and consider as in the proof of Theorem~\ref{Theorem2.4} the process
$
X^\vartheta_{t,u}:= x + \int_t^u \vartheta_r\dd S_r
$,
$t\le u\le T$. (Again, we usually do not explicitly indicate the
dependence of $X^\vartheta$ on the starting time $t$, nor on $x$.) Lemma~\ref{Lem1.5} yields $\vnull= \mnull+ \anull$, and as $\vtwo$ satisfies
the BSDE (\ref{equ3.1}), the same computation as for (\ref{equ2.38}) gives,
with (\ref{equ2.42}), that
\[
(X^\vartheta_u)^2 \vtwo_u - x^2 \vtwo_t
=
m_u - m_t + \int_t^u (\vartheta_r - \gamma_r X^\vartheta_{r-})^2 \N
_r\bigl(\vtwo\bigr)
\dd\MM_r.
\]
Finally, using the product rule, (\ref{equ2.7}), the BSDE (\ref{equ3.2})
for $\vone$, (\ref{equ3.8}) and (\ref{equ2.10}) leads to
\begin{eqnarray*}
d\bigl(\vone X^\vartheta\bigr)
&=&
\vone_- \vartheta\dd S + X^\vartheta_- \dd\vone+ \vartheta\dd\bigl[\vone, S\bigr]
\\
&=&
dm + \vone_- \vartheta\lambda\dd\MM- X^\vartheta_- \gamma
\bigl(\psione+ \lambda
{}^{\p}\vone\bigr) \dd\MM\\
&&{}+ \vartheta\bigl(\psione+ \lambda\Delta\aone\bigr) \dd\MM
\\
&=&
dm + \bigl(\psione+ \lambda{}^{\p}\vone\bigr) (\vartheta
- \gamma X^\vartheta_-) \dd\MM.
\end{eqnarray*}
Using (\ref{equ3.5}) and adding up therefore gives
%
\begin{eqnarray}
\label{equ3.10}
\VH_u(X^\vartheta_u)
&=&
\vnull_u - 2\vone_u X^\vartheta_u + \vtwo_u (X^\vartheta_u)^2
\nonumber\\
&=&
\VH_t(x) + \anull_u - \anull_t\nonumber\\[-8pt]\\[-8pt]
&&{}
- \int_t^u 2 \bigl( \psione_r + \lambda_r {}^{\p}\vone
_r\bigr) (\vartheta_r - \gamma_r
X^\vartheta_{r-}) \dd\MM_r
\nonumber\\
&&{}
+\int_t^u (\vartheta_r - \gamma_r X^\vartheta_{r-})^2 \N_r\bigl(\vtwo\bigr)
\dd\MM_r +
m_u - m_t.
\nonumber
\end{eqnarray}
Now choose $x=0$ and $\vartheta$ of the form $\vartheta= y I_{\rsb t,
\varrho
_t\rsb}$ for some constant $y\in\RR$, where the stopping time
$\varrho
_t>t$ is chosen such that $\vartheta$ is in $\Theta$; this is possible
because $S$ is in $\SS^2_\loc(P)$. Then
$
\vartheta_r = y I_{\{t<r\le\varrho_t\}}
$
and
$
X^\vartheta_{r-} = y (S_{r-}-S_t) I_{\{t<r\le\varrho_t\}}
$,
and plugging this into (\ref{equ3.10}) and collecting terms gives
\begin{eqnarray*}
&&\VH_u(X^\vartheta_u) - \VH_t(0)\\
&&\qquad =
\anull_u - \anull_t
- 2 \int_t^{u\wedge\varrho_t} y \bigl(\psione_r + \lambda_r {}^{\p}\vone_r\bigr) \bigl(
1 - (S_{r-} - S_t) \gamma_r \bigr) \dd\MM_r
\\
&&\qquad\quad{}
+ \int_t^{u\wedge\varrho_t} y^2 \bigl( 1 - (S_{r-} - S_t) \gamma_r
\bigr)^2
\N_r\bigl(\vtwo\bigr) \dd\MM_r + m_u - m_t.
\end{eqnarray*}
By Proposition~\ref{Prop1.1}, this process is always a $P$-submartingale on
$\lsb t,T\rsb$. So if we take a predictable $C\in\A^+_\loc(P)$ with
$\MM
\ll C$ and $\anull\ll C$, we obtain that the process
\begin{eqnarray*}
&&\int_t^{u\wedge\varrho_t} \biggl( \bigl(
y^2 \bigl( 1 - \gamma_r (S_{r-}-S_t) \bigr)^2 \N_r\bigl(\vtwo\bigr)\\
&&\hspace*{12pt}\qquad {}- 2 y \bigl(\psione_r + \lambda_r
{}^{\p}\vone_r\bigr) \bigl( 1 - \gamma_r
(S_{r-}-S_t) \bigr) \bigr) \,{d\MM_r \over dC_r}
+ {d\anull_r \over dC_r} \biggr) \dd C_r
\end{eqnarray*}
for $t\le u\le T$ is, for all $t\in[0,T]$ and $y\in\RR$, an increasing
process. Again arguing as in steps (1) and (2) of the proof of Theorem~\ref{Theorem2.4} and using that $S_{r-}-S_s\to0$ when $s$ increases to $r$
(used for the jumps) or when $r$ decreases to $s$ (used for the
continuous part),
%
%
we get
\begin{eqnarray}
&&y^2 \N\bigl(\vtwo\bigr) \,{d\MM\over dC} - 2 y \bigl(\psione+ \lambda
{}^{\p}\vone\bigr)\,
{d\MM\over dC}
+ {d\anull\over dC}
\ge0\nonumber\\
&&\eqntext{\mbox{for all $y\in\RR$, $P\otimes C$-a.e.}}
\end{eqnarray}
Because $\N(\vtwo) > 0$ by Lemma~\ref{Lem2.3}, we conclude that
%
\begin{equation}\label{equ3.11}
{(\psione+ \lambda{}^{\p}\vone)^2 \over\N
(\vtwo)} {d\MM\over dC}
\le
{d\anull\over dC},
\qquad\mbox{$P\otimes C$-a.e.}
\end{equation}
This implies that
$
\int\{ d\anull- {(\psione+ \lambda\spP{\vone})^2 \over\N
(\vtwo
)} \dd\MM\}
$
is an increasing process, and since $\anull$ is $P$-integrable
because $\vnull$ is a $P$-submartingale by Lemma~\ref{Lem1.5}, we obtain that
\[
E\biggl[ \int_0^T {(\psione_r + \lambda_r {}^{\p}\vone_r)^2 \over\N
_r(\vtwo)} \dd\MM_r \biggr] < \infty.
\]
So if we define
%
\begin{eqnarray}
\label{equ3.12}
\Ynull_t
:\!&=&
E\biggl[ H^2 - \int_t^T {(\psione_r + \lambda_r {}^{\p}\vone_r)^2 \over\N
_r(\vtwo)} \dd\MM_r
\bigg| \F_t \biggr]
\nonumber\\[-8pt]\\[-8pt]
&=&\!:
\Nnull_t + \int_0^t {(\psione_r + \lambda_r {}^{\p}\vone_r)^2 \over\N
_r(\vtwo)} \dd\MM_r,
\nonumber
\end{eqnarray}
then clearly $(\Ynull, \Nnull)$ solves (\ref{equ3.3}), and $\Ynull$ is a
true $P$-submartingale. This shows that there exists a solution to
(\ref{equ3.3})
with (iii), but we do not know yet if $\vnull= \Ynull$.

(3) To finish\vspace*{1pt} the implication ``(a) $\Longrightarrow$
(b),'' we now want to prove that each $X^{\thetast(x,H)}$ satisfies
(\ref{equ3.4}) and that $\vnull= \Ynull$. We again
fix $t$, take $\vartheta\in \Theta$ and do the same calculation as in
(\ref{equ3.10}). Completing the square then gives
%
\begin{eqnarray}
\label{equ3.13}\qquad
\VH_u(X^\vartheta_u)
&=&
\VH_t(x) + m_u - m_t\nonumber\\
&&{}+\int_t^u \biggl( d\anull_r - { (\psione_r + \lambda_r
{}^{\p}\vone_r)^2
\over\N_r(\vtwo) } \dd\MM_r \biggr)
\\
&&{}
+\int_t^u \biggl( (\vartheta_r - \gamma_r X^\vartheta_{r-} ) \sqrt
{ \N_r\bigl(\vtwo
\bigr) }
- { \psione_r + \lambda_r {}^{\p}\vone_r \over
\sqrt{ \N_r(\vtwo) } }
\biggr)^2 \dd\MM_r.
\nonumber
\end{eqnarray}
By Proposition~\ref{Prop1.1}, this process must be a $P$-martingale on
$\lsb t,T\rsb$ if we plug in for $\vartheta$ the optimal $\thetast(x,H)$.
Because both integral terms on the right-hand side are increasing due
to (\ref{equ3.11}), they must then both vanish identically, on $\lsb
t,T\rsb$ for every $t$. This first gives that
%
\begin{equation}\label{equ3.14}
\anull= \int{(\psione+ \lambda{}^{\p}\vone)^2
\over\N(\vtwo)} \dd\MM,
\end{equation}
and as $\vnull= \mnull+ \anull$ is a $P$-submartingale, comparing
(\ref{equ3.12}) and (\ref{equ3.14}) yields $\mnull= \Nnull$, hence $\vnull= \Ynull$, and
so $(\vnull, \mnull)$ solves the BSDE (\ref{equ3.3}) and also is the unique
solution satisfying (iii).
Second, we obtain for the optimal strategy $\thetast= \thetast(x,H)$ that
\[
\thetast_u
=
{\psione_u + \lambda_u {}^{\p}\vone_u \over\N
_u(\vtwo)} + \gamma_u
X^{\thetast}_{u-},
\]
which is (\ref{equ3.6}) in view of the definition (\ref{equ2.42}) of $\gamma
$; recall that $(\vtwo, \psitwo, \Ltwo)$ solves (\ref{equ2.18}).
Integrating with respect to $S$ shows that $X^{\thetast}$ satisfies the
SDE (\ref{equ3.4}) on $\lsb t,T\rsb$, and since $\thetast$ is in $\Theta
$, the unique solution of (\ref{equ3.4}) is in $\SS^2(P)$. So we have
now proved that (a) implies (b), and also that we then have (\ref{equ3.5})
and (\ref{equ3.6}).

(4) Conversely, let us start with (b); then we have to prove the
existence of an optimal $\thetast(x,H)$. Fix $t$, set $ W_u(x):=
\Ynull_u - 2\Yone_u x + \Ytwo_u x^2 $ for $t\le u\le T$ and use
(\ref{equ2.22}) and the BSDEs (\ref{equ3.1})--(\ref{equ3.3}) for $\Ytwo, \Yone, \Ynull$ to compute
as for (\ref{equ3.10}) and (\ref{equ3.13}) that for any $\vartheta\in\Theta$,
%
\begin{eqnarray}\label{equ3.15}\qquad
W_u(X^\vartheta_u)
&=&
W_t(x) + m_u - m_t\nonumber\\[-8pt]\\[-8pt]
&&{}+\int_t^u \biggl( (\vartheta_r - \gamma_r X^\vartheta_{r-} )
\sqrt{ \N_r\bigl(\Ytwo
\bigr) }
- { \psione_r + \lambda_r {}^{\p}\Yone_r \over
\sqrt{ \N_r(\Ytwo) } }
\biggr)^2 \dd\MM_r
\nonumber
\end{eqnarray}
for $t\le u\le T$. So $W(X^\vartheta)$ is a local $P$-submartingale on
$\lsb t,T\rsb$; but we also know from (b) that $\Ynull$
is of class (D), $\Ytwo$ is bounded, and $|\Yone|^2$ is of class (D).
Since $X^\vartheta$ is in $\SS^2(P)$ for every $\vartheta\in\Theta
$, we see
that $W(X^\vartheta)$ is thus of class (D), hence a true
$P$-submartingale, and so
\[
W_t(x)
\le
E[ W_T(X^\vartheta_T) | \F_t ]
=
E\biggl[ \biggl( H - x - \int_t^T \vartheta_r\dd S_r \biggr)^2
\bigg| \F_t
\biggr]
\]
for any $\vartheta\in\Theta$. This yields $W_t(x) \le\VH_t(x)$ by (\ref{equ1.2}).
Conversely, if we take the solution $X^{(t)}$ of (\ref{equ3.4}) and define
\[
\thetat^{(t)}
:=
{\psione+ \lambda{}^{\p}\Yone\over\N(\Ytwo
)} - {\psitwo+ \lambda
{}^{\p}\Ytwo\over\N(\Ytwo)} X^{(t)}_-,
\]
then\vspace*{1pt} integrating with respect to $S$ shows that $
X^{\thetat^{(t)}} = x + \int_t^{\fatdot} \thetat^{(t)}_r\dd S_r $
equals $X^{(t)}$, since both satisfy (\ref{equ3.4}), and is in $\SS
^2(P)$ due to (b) so that $\thetat^{(t)}$ is in $\Theta$. Moreover,
plugging in $\thetat^{(t)}$ for $\vartheta$ shows, similar to the argument for
(\ref{equ3.15}), that $W(X^{\thetat^{(t)}})$ is a (true) $P$-martingale
on $\lsb t,T\rsb $. This implies that
\[
W_t(x)
=
E\biggl[ \biggl( H - x - \int_t^T \thetat^{(t)}_r\dd S_r \biggr)^2
\bigg|
\F_t \biggr]
\ge
\VH_t(x),
\]
and so\vspace*{1pt} we conclude that $W_t(x) = \VH_t(x)$ and that $\thetat^{(t)}$ is
optimal for (\ref{equ1.2}), giving existence of
$\thetast(x,H):= \thetat^{(t)}$. This proves that (b) implies (a) and
that we then also have $W(x) = \VH(x)$ for all $x$, hence $Y^{(i)} =
\vi
$ for $i=0,1,2$. This ends the proof of (1).

(5) Finally, the assertion of part (2) follows, similarly to Theorem~\ref{Theorem2.4},
from the proof of part (1); we only need to notice again that $L^2(\F
_t,P) + G_{t,T}(\Theta)$ is closed in $L^2(P)$ for every $t$.
\end{pf*}

%
%
\section{Alternative versions for the BSDEs}\label{sec4}

In this section, we give equivalent alternative versions for the BSDEs
obtained in Sections~\ref{sec2} and~\ref{sec3}. One reason is that in some
models, these versions are more convenient to work with; a second is
that it allows us to discuss how our results relate to existing literature.

For reasons of space, we only look at (\ref{equ2.18}) or (\ref{equ3.1}) in
detail; this is the most complicated equation. \textit{Throughout this
section, we assume as in Theorem~\ref{Theorem2.4} that $S\in\SS^2_\loc(P)$ and
$\PPes\ne\varnothing$.} For convenience, we recall that (\ref{equ2.18}) reads
%
\begin{equation}\label{equ4.1}\quad
Y_t = Y_0 + \int_0^t {( \psi_s + \lambda_s {}^{\p}Y_s )^2 \over\NNY{s} }
\dd\MM_s +
\int_0^t \psi_s\dd M_s + L_t,\qquad
Y_T=1,
\end{equation}
where
$
{\mathcal N}(Y) = {}^{\p}Y (1+\lambda^2\Delta\MM) + \func(Y)
$
and
$
\func(Y) = {d\sdpP{[ \NY, [S] ]} \over d\MM},
$
as in (\ref{equ2.12}) and (\ref{equ2.11}). A solution of (\ref{equ4.1}) is a
priori a tuple $(Y,\psi,L)$ with $L\in\M_{0, \loc}(P)$ strongly
$P$-orthogonal to $M$, $\psi\in L^1_\loc(M)$, and $Y$ a $P$-special
semimartingale such that $[ \NY, [S] ]\in\A_\loc(P)$. In view of
Theorem~\ref{Theorem2.4} (where $Y$ is bounded), we restrict ourselves to
solutions with $\psi\in L^2_\loc(M)$ and $L\in\M^2_{0,\loc}(P)$. For
better comparison with (\ref{equ3.1}), we really ought to write a
superscript ${}^{(2)}$ for $Y, \psi, L$, but we omit this to alleviate
the notation.

%
%
\subsection{Working with $M^d$}\label{sec4.1}
%
%

%
The BSDE (\ref{equ4.1}) is written with the local $P$-martingale $M$ from
the canonical decomposition
$
S = S_0 + M + A = S_0 + M + \int\lambda\dd\MM
$
of $S$. In simple models with jumps, it is useful to split $M = M^c +
M^d$ into its continuous and purely discontinuous local martingale
parts $M^c$ and $M^d$, respectively. Then
$
\MM= \MMc+ \MMd
$,
and we define the predictable processes
\[
\deltac:= {d\MMc\over d\MM},\qquad
\deltad:= {d\MMd\over d\MM} = 1 - \deltac.
\]

We now consider the backward equation
%
\begin{eqnarray}
\label{equ4.2}
Y_t
&=&
Y_0 + \int_0^t {( \psic_s \deltac_s + \psid_s (1-\deltac_s) +
\lambda_s {}^{\p}Y_s )^2 \over{}^{\p}Y_s (1+\lambda_s^2\Delta\MM_s) +
\func_s(Y)} \dd\MM_s\nonumber\\
&&{}
+ \int_0^t \psic_s \dd\Mc_s + \int_0^t \psid_s \dd\Md_s + \Ls_t,
\\
Y_T&=&1.
\nonumber
\end{eqnarray}
A \textit{solution} of (\ref{equ4.2}) is a priori a tuple $(Y, \psic, \psid
, \Ls)$ with $\Ls\in\M_{0,\loc}(P)$ strongly $P$-orthogonal to
both $\Mc
$ and $\Md$, $\psic\in L^2_\loc(\Mc)$, $\psid\in L^1_\loc(\Md)$ and
$Y$ a \mbox{$P$-special} semimartingale with $[ \NY, [S] ] \in\A_\loc(P)$.
As for (\ref{equ4.1}), we restrict our attention to solutions with $\psid
\in L^2_\loc(\Md)$ and $\Ls\in\M^2_{0,\loc}(P)$.
%
\begin{proposition}\label{Prop4.1}
The BSDEs (\ref{equ4.1}) and (\ref{equ4.2}) are equivalent. More precisely, $(Y,\psi,L)$
with $\psi\in L^2_\loc
(M)$ and $L\in\M^2_{0,\loc}(P)$ solves (\ref{equ4.1}) if and only if
$(Y,\psic,\psid,\Ls)$ with $\psic\in L^2_\loc(\Mc)$, $\psid\in
L^2_\loc
(\Md)$ and $\Ls\in\M^2_{0,\loc}(P)$ solves~(\ref{equ4.2}), where the
tuples are related by
%
\begin{equation}\label{equ4.3}
\psi\SIdot M + L = \psic\SIdot\Mc+ \psid\SIdot\Md+ \Ls.
\end{equation}
\end{proposition}
\begin{pf}
If $(Y,\psi,L)$ solves (\ref{equ4.1}), we use the
Galtchouk--Kunita--Watanabe decomposition of $\psi\SIdot M + L$ with
respect to $\Mc$ and $\Md$ to obtain (\ref{equ4.3}) and define $\psic,
\psid, \Ls$; so $\Ls$ is strongly $P$-orthogonal to both $\Mc$ and
$\Md
$, and taking the covariation with $M$ and using $\langle L, M\rangle
\equiv0$ gives $\psi= \psic\deltac+ \psid\deltad$. Plugging this
and (\ref{equ4.3}) into (\ref{equ4.1}) shows directly that $(Y,\psic,\psid
,\Ls)$ solves~(\ref{equ4.2}).

Conversely, if $(Y,\psic,\psid,\Ls)$ solves (\ref{equ4.2}), we define
\[
\psi:= \psic\deltac+ \psid(1-\deltac)\in L^2_\loc(M)
\]
and $L:= \psic \SIdot\Mc+ \psid\SIdot\Md+ \Ls- \psi\SIdot M
\in\M^2_{0,\loc}(P)$. Then plugging into (\ref{equ4.2}) directly shows
that $(Y,\psi,L)$ satisfies (\ref{equ4.1}), and since $\langle
L,M\rangle\equiv0$ due to the definitions above, $L$ is also strongly
$P$-orthogonal to $M$. So $(Y,\psi,L)$ solves~(\ref{equ4.1}).
\end{pf}

Equation (\ref{equ4.2}) is particularly convenient for models with simple
jumps, as illustrated by:
%
\begin{example}\label{Exa4.2}
Consider the \textit{jump-diffusion model}
\[
dS_t = S_{t-} (\mu_t\dd t + \sigma_t \dd W_t + \eta_t \dd n_t),\qquad
S_0>0,
\]
where $W$ is a Brownian motion, and $n_t = N_t - \alpha t$, $0\le t\le
T$, is the compensated martingale of a simple Poisson process with
intensity $\alpha>0$. The predictable processes $\mu, \sigma, \eta$
satisfy $\sigma\ne0$ and suitable integrability conditions, and we
assume that $\eta>-1$ to ensure that $S>0$. Then\vspace*{-1pt} we have
$
d\Mc_t = S_{t-} \sigma_t \dd W_t
$,
$
d\Md_t = S_{t-} \eta_t \dd n_t
$,
$
d\MM_t = S_{t-}^2 (\sigma_t^2 + \alpha\eta_t^2) \dd t
$,
$
\lambda_t = {\mu_t \over S_{t-} (\sigma_t^2 + \alpha\eta_t^2)}
$\vadjust{\goodbreak}
and
$
\deltac_t = {\sigma_t^2 \over\sigma_t^2 + \alpha\eta_t^2}
$.
Because $\MM$ is continuous, so is $B^Y$ due to (\ref{equ4.2}); hence
${}^{\p}Y = Y_-$ by (\ref{equ2.10}). Moreover, using $[n] = N$ gives
\begin{eqnarray*}
\dpP{[ \NY, [S] ]}_t
&=&
\dpP{\bigl[ \psic\SIdot\Mc+ \psid\SIdot\Md+ \Ls, [\Md] \bigr]}_t
\\[-2pt]
&=&
\dpP{\bigl[ \psid\SIdot\Md+ \Ls, (S_- \eta)^2\SIdot[n] \bigr]}_t
\\[-2pt]
&=&
(S_-^3 \psid\eta^3)\SIdot\dpP{N}_t
\\[-2pt]
&=&
(S_-^3 \psid\eta^3 \alpha)\SIdot t
\end{eqnarray*}
so that
$
\func_t(Y) = {\alpha\eta_t^3 \psid_t S_{t-} \over\sigma_t^2 +
\alpha
\eta_t^2}
$.
Using the notation $\widetilde\psi^c= \psic S_- \sigma$, $\widetilde\psi^d=
\psid S_- \eta$ and plugging in then allows us to rewrite the BSDE
(\ref{equ4.2}) after simple calculations as
\begin{eqnarray*}
Y_t &=& Y_0 + \int_0^t {(\widetilde\psi^c_s \sigma_s +
\alpha\widetilde\psi^d_s \eta_s + \mu_s Y_{s-} )^2 \over Y_{s-}
(\sigma_s^2 + \alpha\eta_s^2) +\alpha\widetilde\psi^d_s \eta
_s^2 } \dd s
+ \int_0^t \widetilde\psi^c_s\dd W_s + \int_0^t \widetilde\psi^d_s
\dd n_s +\Ls_t,\\[-2pt]
Y_T &=& 1.
\end{eqnarray*}
It depends on the choice of the filtration $\FF$ whether we can have a
nontrivial $\Ls\in\M^2_{0,\loc}(P)$ strongly $P$-orthogonal to both
$\Mc
$ and $\Md$, or $W$ and $n$. If $\FF$ is generated by $W$ and $N$, then
$\Ls\equiv0$ automatically by the martingale representation theorem in
$\FF^{W,N}$.\vspace*{-2pt}
\end{example}

%
\subsection{Using random measures}\label{sec4.2}
%
%

%
For models with more general jumps, the version (\ref{equ4.2}) of the
basic BSDE (\ref{equ4.1}) is less useful because one cannot easily
express $\func(Y)$ in terms of integrands like in the preceding
example. We therefore use semimartingale characteristics and, in
particular, work with the jump measure of $S$. For the required
notation and results, we refer to Chapter II of \citet{JacShi03}. We take $E=\RR$
there so that $\widetilde\Omega= \Omega\times[0,T]\times\RR$ with the
$\sigma$-field $\widetilde\P= \P\otimes\B(\RR)$, where $\P$ is the
predictable $\sigma$-field on $\Omega\times[0,T]$.

Denote by $\mu^S$ the random measure associated with the jumps of $S$
and by $\nu$ its $P$-compensator. Using Proposition II.2.9 of
\citet{JacShi03},
we have
\[
\nu(\omega, dt, dx) = F_t(\omega, dx) \dd B_t(\omega)
\]
for a predictable increasing $B$ null at 0. Moreover, (\ref{equ2.7}) gives
$\Delta S = \Delta M + \lambda\Delta\MM$ and $(x^2\wedge1)*\mu^S
\ll
[M]+\MM$, and combining this with the construction of $B$ in \citet{JacShi03} and
(\ref{equ2.6}), we see that $B\ll\MM$. We introduce the predictable processes
\[
b:= {dB \over d\MM},\qquad
\deltac:= {d\MMc\over d\MM}
\]
and note that
$
[\Md] = \sum(\Delta M)^2 = (x-\lambda\Delta\MM)^2 * \mu^S
$
implies that
\[
\MMd= (x-\lambda\Delta\MM)^2 *\nu= \biggl( \int(x-\lambda\Delta
\MM)^2
F(dx) \biggr)\SIdot B,\vadjust{\goodbreak}
\]
so that $\MM= \MMc+ \MMd$ can be reformulated as
%
\begin{equation}\label{equ4.4}
\deltac_t + b_t \int(x-\lambda_t\Delta\MM_t)^2 F_t(dx) = 1,
\qquad\mbox{$P\otimes\MM$-a.e.}
\end{equation}

With the notation
$
\Wh_t = \int_{\RR} W_t(x) \nu(\{t\}, dx)
$,
we now consider the backward equation
%
\begin{eqnarray}\quad
\label{equ4.5}
Y_t
&=&
Y_0 + \int_0^t {( \varphi_s\deltac_s + b_s \int x (W_s(x)-\Wh_s)
F_s(dx) + \lambda_s {}^{\p}Y_s )^2
\over{}^{\p}Y_s \deltac_s + b_s \int x^2 (
{}^{\p}Y_s + W_s(x)-\Wh_s
) F_s(dx)} \dd\MM_s
\nonumber\\[-8pt]\\[-8pt]
&&{}
+ \int_0^t \varphi_s \dd\Mc_s + W*(\mu^S-\nu)_t + \Ls_t,\qquad
Y_T = 1.
\nonumber
\end{eqnarray}
A \textit{solution} of (\ref{equ4.5}) is a priori a tuple $(Y,\varphi,W,\Ls)$
such that $\varphi\in L^2_\loc(\Mc)$, \mbox{$W\in\G^1_\loc(\mu^S)$} [see
(3.62) in
\citet{Jac79}], $\Ls\in\M^2_{0,\loc}(P)$ strongly $P$-ortho\-gonal to
$\Mc$ and
to the space of stochastic integrals
$
\{ \bar W*(\mu^S-\nu) | \bar W\in\G^2_\loc(\mu^S) \}
$,
and $Y$ a $P$-special semimartingale with $[ \NY, [S] ]\in\A_\loc
(P)$. As before for (\ref{equ4.1}) and (\ref{equ4.2}), we restrict our
attention to solutions with $W\in\G^2_\loc(\mu^S)$ and $\Ls\in\M
^2_{0,\loc}(P)$.\vspace*{1pt}

In view of the next result, (\ref{equ4.5}) seems the natural form of the
BSDE (\ref{equ4.1}) or (\ref{equ2.18}) in the general case, because its
generator is expressed in terms of integrands. Nevertheless, as seen in
Section~\ref{sec2}, the form (\ref{equ2.18}) is more convenient for proving
results via stochastic calculus.
%
\begin{proposition}\label{Prop4.3}
The BSDEs (\ref{equ4.1}) and (\ref{equ4.5}) are equivalent. More precisely, $(Y,\psi,L)$
with $\psi\in L^2_\loc
(M)$, and $L\in\M^2_{0,\loc}(P)$ solves (\ref{equ4.1}) if and only if
$(Y,\varphi,W,\Ls)$ with $\varphi\in L^2_\loc(\Mc)$, $W\in\G
^2_\loc(\mu^S)$
and $\Ls\in\M^2_{0,\loc}(P)$ solves (\ref{equ4.5}), where the tuples are
related by the equation
\[
\psi\SIdot M + L = \varphi\SIdot\Mc+ W*(\mu^S-\nu) + \Ls.
\]
\end{proposition}
\begin{pf}
If $(Y,\psi,L)$ solves (\ref{equ4.1}), we take its martingale part
$\psi\SIdot M + L$ and represent this as
%
\begin{equation}\label{equ4.6}
\psi\SIdot M + L = \varphi\SIdot\Mc+ W*(\mu^S-\nu) + U*\mu^S + \Lt
\end{equation}
with $\varphi\in L^2_\loc(\Mc)$, $W\in\G^2_\loc(\mu^S)$, $U\in
\H^2_\loc(\mu ^S)$ [see\vspace*{1pt} \citet{Jac79}, Section~3.3b,
pages 101 and 102] and $\Lt\in\M ^2_{0,\loc}(P)$ with $[ \Lt,S ]
\equiv0$. This is the so-called \textit{Jacod decomposition}; see
\citet{Jac79}, Theorem 3.75, or Theorem~2.4 in
\citet{ChoSch11} for a more detailed exposition.

We next express $\func(Y)$ in terms of $W$ and $\nu$. Using (\ref{equ4.1}) and
(\ref{equ4.6}) yields
%
\begin{equation}\label{equ4.7}
\Delta\NY_t = W_t(\Delta S_t) I_{\{\Delta S_t\ne0\}} - \Wh_t +
U_t(\Delta S_t) I_{\{\Delta S_t\ne0\}} + \Delta\Lt_t.
\end{equation}
Moreover,
$
\sum\Delta\Lt(\Delta S)^2 = \Delta S\SIdot[ \Lt,S ] \equiv0
$
so that we get
\[
[ \NY, [S] ]
=
\sum\Delta\NY(\Delta S)^2
=
\bigl( x^2 \bigl(W(x) - \Wh\bigr) \bigr) *\mu^S + ( x^2 U(x) )*\mu^S.
\]
Because $[ \NY, [S] ]$ is in $\A_\loc(P)$, this implies that $x^2
U(x)$ is in $\H^1_\loc(\mu^S)$ so that $( x^2 U(x))*\mu^S$ is a local
$P$-martingale by \citet{Jac79}, (3.73). Hence we obtain
\begin{eqnarray*}
\dpP{[ \NY, [S] ]}
&=&
\dpP{\bigl( \bigl( x^2 \bigl(W(x) - \Wh\bigr) \bigr) *\mu^S \bigr)}
=
\bigl( x^2 \bigl(W(x) - \Wh\bigr) \bigr)*\nu\\
&=&
\biggl( \int x^2 \bigl(W(x) - \Wh\bigr) F(dx) \biggr)\SIdot B,
\end{eqnarray*}
and so
$
\func_t(Y) = b_t \int x^2 (W_t(x) - \Wh_t) F_t(dx)
$.
Moreover,
\[
[S] = [S]^c + \sum(\Delta S)^2 = \MMc+ x^2 * \mu^S
\]
gives
$
\dpP{[S]} = \MMc+ x^2 * \nu= ( \deltac+ \int x^2 F(dx) b) \SIdot
\MM
$
so that comparing with (\ref{equ2.8}) yields that
$
1 + \lambda^2 \Delta\MM= \deltac+ b \int x^2 F(dx)
$
and hence
%
\begin{eqnarray}\label{equ4.8}
\NNY{t}
&=&
{}^{\p}Y_t (1+\lambda_t^2 \Delta\MM_t) + \func
_t(Y)\nonumber\\[-8pt]\\[-8pt]
&=&
{}^{\p}Y_t \deltac_t + b_t \int x^2 \bigl(
{}^{\p}Y_t + W_t(x)-\Wh_t \bigr) F_t(dx).
\nonumber
\end{eqnarray}
If we now define $\Ls:= U*\mu^S + \Lt$, then (\ref{equ4.6}) gives
%
\begin{equation}\label{equ4.9}
\psi\SIdot M + L = \varphi\SIdot\Mc+ W*(\mu^S-\nu) + \Ls.
\end{equation}
But $ [\Ls,M] = [\Ls, S] - [\Ls, \lambda\SIdot\MM] =( x U(x)) *\mu^S +
[ \Lt ,S ] - [\Ls, \lambda\SIdot\MM] $ is a local $P$-martingale by
Yoeurp's lemma, and a similar argument as just above, using now that
$U\in\H^2_\loc(\mu^S)$; so $\langle\Ls, M\rangle\equiv0$ and $\Ls$ is
strongly $P$-orthogonal to~$\Mc$. Moreover, we have for all
\mbox{$\bar W\in\G^2_\loc(\mu^S)$} that $ [ \Lt, \bar W*(\mu^S-\nu) ] =
0 $ since $[ \Lt,S ]\equiv0$, and so\vspace*{1pt} $ \langle\Ls, \bar W*(\mu^S-\nu)
\rangle= \langle U*\mu^S, \bar W*(\mu ^S-\nu) \rangle\equiv0 $ for
all
$\bar W\in\G^2_\loc(\mu^S)$ by \citet{Jac79}, Exercice 3.23.
Finally, (\ref{equ2.7}) and Yoeurp's lemma yield
%
\begin{eqnarray}
\label{equ4.10}
\langle W*(\mu^S-\nu), M\rangle
&=&
\dpP{[ W*(\mu^S-\nu), S - \lambda\SIdot\MM]}
\nonumber\\
&=&
\dpP{[ W*(\mu^S-\nu), S ]}
\nonumber\\[-8pt]\\[-8pt]
&=&
\dpP{\bigl( \bigl( x \bigl(W(x) - \Wh\bigr) \bigr) * \mu^S \bigr)}
\nonumber\\
&=&
\bigl( x \bigl(W(x) - \Wh\bigr) \bigr) * \nu.
\nonumber
\end{eqnarray}
Taking in (\ref{equ4.9}) the covariation with $M$ and using also $\langle
L, M\rangle\equiv0 \equiv\langle\Ls, M\rangle$ yields
\[
\psi\SIdot\MM
=
\biggl( \varphi\deltac+ \biggl( \int x \bigl(W(x)-\Wh\bigr) F(dx) \biggr) b
\biggr)\SIdot
\MM
\]
so that we get
%
\begin{equation}\label{equ4.11}
\psi_t = \varphi_t \deltac_t + b_t \int x \bigl(W_t(x) - \Wh_t\bigr) F_t(dx),
\qquad\mbox{$P\otimes\MM$-a.e.}
\end{equation}
Plugging (\ref{equ4.11}) and (\ref{equ4.8}) into (\ref{equ4.1}) and using
(\ref{equ4.9}), we see that $(Y,\varphi,W,\Ls)$ solves (\ref{equ4.5}).

Conversely, if $(Y,\varphi,W,\Ls)$ solves (\ref{equ4.5}), then we define
$\psi$ by (\ref{equ4.11}) and
\[
L:= \varphi\SIdot\Mc- \psi\SIdot M + W*(\mu^S-\nu) + \Ls.
\]
Then\vspace*{1pt} $\psi\in L^2_\loc(M)$, due to (\ref{equ4.4}) and because $W\in\G
^2_\loc(\mu^S)$, and so $L\in\M^2_{0,\loc}(P)$. Moreover, equation
(\ref{equ4.10}), the
definitions of $L$ and $\psi$ via (\ref{equ4.11}) and the definitions of
$\deltac$ and $b$ yield
\[
\langle L, M\rangle
=
\langle\Ls, M\rangle
=
\langle\Ls, \Mc+\Md\rangle
=
\langle\Ls, \Mc\rangle+ \langle\Ls, x*(\mu^S-\nu)\rangle
\equiv
0
\]
by the orthogonality properties of $\Ls$, so that $L$ is strongly
$P$-orthogonal to $M$. Finally, the Jacod\vspace*{1pt} decomposition applied to $\Ls
$ implies that the latter must have the form $\Ls= U*\mu^S + \Lt$ due
to its orthogonality properties. But then we obtain from (\ref{equ4.5})
again (\ref{equ4.7}), hence also (\ref{equ4.8}), and then plugging in shows
that $(Y,\psi,L)$ solves (\ref{equ4.1}). This completes the proof.
\end{pf}

Just for completeness, but without any details, we give here the
equivalent versions of the BSDEs (\ref{equ3.2}) and (\ref{equ3.3}) for
$\vone$ and $\vnull$. They are
\begin{eqnarray*}
d\Yone_t
&=&
{( \phione_t\deltac_t + b_t \int x( \Wone_t(x) - \widehat
{W}^{(1)}_t
) F_t(dx) + \lambda_t {}^{\p}\Yone_t )
\over
{}^{\p}\Ytwo_t \deltac_t + b_t \int x^2 (
{}^{\p}\Ytwo_t + \Wtwo_t(x) -
\widehat{W}^{(2)}_t ) F_t(dx)}
\\
&&{}
\times\biggl( \phitwo_t\deltac_t + b_t \int x \bigl(\Wtwo_t(x) - \widehat
{W}^{(2)}_t\bigr) F_t(dx) + \lambda_t {}^{\p}\Ytwo_t
\biggr) \dd\MM_t
\\
&&{}
+ \phione_t\dd\Mc_t + d\bigl( \Wone*(\mu^S-\nu) \bigr)_t + dL^{(1),
\prime}_t,\qquad
\Yone_T = H,
\end{eqnarray*}
and
\begin{eqnarray*}
d\Ynull_t
&=&
{( \phione_t\deltac_t + b_t \int x (\Wone_t(x) - \widehat
{W}^{(1)}_t)
F_t(dx) + \lambda_t {}^{\p}\Yone_t )^2
\over
{}^{\p}\Ytwo_t \deltac_t + b_t \int x^2 (
{}^{\p}\Ytwo_t+ \Wtwo_t(x) -
\widehat{W}^{(2)}_t ) F_t(dx)} \dd\MM_t + d\Nnull_t,
\\
\Ynull_T &=& H^2.
\end{eqnarray*}
Finally, the recursive representation for the optimal strategy in
(\ref{equ3.6}) takes the form
\begin{eqnarray*}
\thetasn_t
&=&
{ \phione_t\deltac_t + b_t \int x (\Wone_t(x) - \widehat{W}^{(1)}_t)
F_t(dx) + \lambda_t {}^{\p}\Yone_t
\over
{}^{\p}\Ytwo_t \deltac_t + b_t \int x^2 (
{}^{\p}\Ytwo_t + \Wtwo_t(x) -
\widehat{W}^{(2)}_t ) F_t(dx)}
\\
&&{}
- { \phitwo_t\deltac_t + b_t \int x (\Wtwo_t(x) - \widehat{W}^{(2)}_t)
F_t(dx) + \lambda_t {}^{\p}\Ytwo_t
\over
{}^{\p}\Ytwo_t \deltac_t + b_t \int x^2 (
{}^{\p}\Ytwo_t + \Wtwo_t(x) -
\widehat{W}^{(2)}_t ) F_t(dx)} X^{\thetasn}_{t-}.
\end{eqnarray*}
Of course, this can equivalently be rewritten as a linear SDE for
$X^{\thetasn}$ as in (\ref{equ3.4}), simply by integrating with respect
to $S$.

%
%
\subsection{Further comments}\label{sec4.3}
%
%

%
At this point, it seems appropriate to comment on related work in the
literature, where we restrict ourselves to papers that have used BSDE
techniques in the context of mean-variance hedging. While extending
work by many authors done for an It\^o process setting in a Brownian
filtration, the results in Mania and Tevzadze (\citeyear{ManTev03N1,ManTev03N2}) and
\citet{BobSch04} still all assume that $S$ is continuous. At the other end of the
scale, \citet{CerKal07} have a general $S\in\SS^2_\loc
(P)$, with
$\PPes\ne\varnothing$; but their methods do not exploit stochastic
control ideas and results at all, and BSDEs appear only very
tangentially in their equations (3.32) and (3.37). As a matter of fact,
their opportunity process $L$ equals our coefficient $\vtwo$, and so
their equation (3.37), which gives a BSDE for $L$, should coincide with
our equation (\ref{equ4.5}). However, \citet{CerKal07} give no
proof for (3.37) and even
remark that ``it is not obvious whether this representation is of
any use.'' Moreover, a closer examination shows that (3.37) is not
entirely correct; it seems that they dropped the jumps of the FV part
of $L$ somewhere, which explains why their equation has $L_-$ instead
of (the correct term)~${}^{\p}L$.

The paper closest to our work is probably \citet{KohXioYe10}. They first study the
variance-optimal martingale measure as in \citet{ManTev03N2}
via the problem
dual to mean-variance hedging and obtain a BSDE that describes $\Vt=
1/\Vn(1) = 1/\vtwo$; see our Proposition~\ref{Prop2.2}. For mean-variance
hedging itself, they subsequently describe the optimal strategy in
feedback form with the help of a process (called $h$) for which they
give a BSDE. Their assumptions are considerably\vspace*{1pt} more restrictive than
ours because, in addition to $S\in\SS^2_\loc(P)$ and $\PPes\ne
\varnothing
$, they also suppose that $S$ is quasi-left-continuous; and for the
results on mean-variance hedging, they additionally even assume that
$\M
^d_\loc(P)$ is generated by integrals of $\mu^S-\nu$ [and also that the
VOMM exists and satisfies the reverse H\"older inequality $R_2(P)$ and
a certain jump condition]. We found it hard to see exactly why this
restrictive condition on $\M^d_\loc(P)$ is needed; the proof in
\citet{KohXioYe10}
for their verification result is rather computational and does not
explain where the rather technical BSDEs come from.

Finally, a similar (subjective) comment as the last one also applies to
\citet{Lim05}. The problem studied there is mean-variance hedging (not the
VOMM), and the process $S$ is a multivariate version of the simple
jump-diffusion model in Example~\ref{Exa4.2}, with a $d$-dimensional
Brownian motion $W$ and an $m$-variate Poisson process $N$. The
filtration used for strategies $\vartheta$ and payoffs $H$ is
generated by
$W$ and $N$; but all model coefficients (including the intensity of
$N$) are assumed to be $\FF^W$-predictable. Technically speaking, this
condition serves to simplify Lim's equation (\ref{equ3.1}), which corresponds
to our equation from Example~\ref{Exa4.2} for $Y$ without the jump term.
It would be interesting to see also at the conceptual level why the
assumption is needed.
\begin{Rem*}
As already pointed out before Theorem~\ref{Theorem3.1}, the BSDE system
(\ref{equ3.1})--(\ref{equ3.3}) is less complicated than it looks. It is
only weakly coupled, meaning that one can solve (\ref{equ3.3}) (even
directly) once one has the solutions of (\ref{equ3.1})
and~(\ref{equ3.2}), and that (\ref{equ3.2}) is linear and hence also
readily solved once one has the solution of~(\ref{equ3.1}). In general,
however, (\ref{equ3.1}) has a very complicated driver, and it seems a
genuine challenge for abstract BSDE theory to prove existence of a
solution directly via BSDE techniques. We do not do that (and do not
need to) since we only use the BSDEs to describe optimal strategies;
existence of the latter (and hence existence of solutions to the BSDEs)
is proved directly via other arguments.

In the special case where the filtration $\FF$ is continuous, the
complicated equation (\ref{equ3.1}) or (\ref{equ2.18}) can be reduced to a
classical quadratic BSDE, as follows. First of all, as already pointed
out before Lemma~\ref{Lem2.3}, the operation ${\mathcal N}(Y)$ in (\ref{equ2.12})
reduces to ${\mathcal N}(Y) = Y$, at least in the context of (\ref{equ2.18}). So
(\ref{equ2.18}) becomes
%
\begin{equation}\label{equ4.12}
dY_t = {(\psi_t + \lambda_t Y_t)^2\over Y_t} \dd\MM_t + \psi_t \dd M_t
+ dL_t,\qquad Y_T=1,
\end{equation}
and we know from Lemma~\ref{Lem2.1} that the solution $\q= \Vn(1)$ is
strictly positive. If we introduce $y:=\log Y$, apply It\^o's
formula and define $\varphi:= \psi/Y$, $\ell:= \int(1/Y)\dd L$,
then it
is straightforward to verify that (\ref{equ4.12}) can be rewritten as
\[
dy_t = \varphi_t \dd M_t + \bigl( (\varphi_t+\lambda_t)^2 - \tfrac{1}{
2} \varphi_t^2
\bigr) \dd\MM_t + d\ell_t - \tfrac{1}{2}\dd\langle\ell\rangle_t,\qquad
y_T=0.
\]
This can then be tackled by standard BSDE methods, if desired.
\end{Rem*}

%
%
\section{Examples}\label{sec5}

In this section, we present some simple examples and special cases to
illustrate our results. We keep this deliberately short in view of the
total length of the paper. \textit{Throughout this section}, \textit{we assume that
$S\in\SS^2_\loc(P)$ and $\PPes\ne\varnothing$.}

Recall the $P$-canonical decomposition $S = S_0 + M + \int\lambda\dd\MM
$ of our price process. Because $\lambda\in L^2_\loc(M)$, the process
$\Zh:= \E(-\lambda\SIdot M)$ is in $\M^2_\loc(P)$ with $\Zh_0=1$.
Moreover, it is easy to check that $\Zh S$ is a local $P$-martingale so
that $\Zh$ is a so-called signed local martingale density for $S$. If
$\Zh$ is a true $P$-martingale and in $\M^2(P)$, then $\Qh$ with $d\Qh
:= \Zh_T\dd P$ is in $\PPss$ and called the \textit{minimal signed}
(\textit{local}) \textit{martingale measure} for $S$; if even $\Zh>0$
so that $\Qh$ is in $\PPes$, then $\Qh$ is the \textit{minimal
martingale measure} (\textit{MMM}) for $S$.\vadjust{\goodbreak}

The MMM is very convenient because its density process $\Zh$ can be
read off explicitly from $S$. On the other hand, the important quantity
for mean-variance hedging is the variance-optimal martingale measure
(VOMM) $\Qt$. By Proposition~\ref{Prop2.6}, we could construct a solution
to the BSDE (\ref{equ2.18}) from $\Qt$ by
\[
\Vn_t(1) = \q_t = \vtwo_t = 1/\Vt_t = {(Z^{\Qt}_t)^2 \over E[
(Z^{\Qt}_T)^2 | \F_t ]},
\qquad 0\le t\le T,
\]
but the density process $Z^{\Qt}$ is usually difficult to find. An
exception is the case when $\Qt= \Qh$, since then $Z^{\Qt} = \Zh=
\E(
-\lambda\SIdot M)$ and the above formula allows us to find an explicit
expression for $\vtwo$. To make this approach work, we need conditions
when $\Qt$ and $\Qh$ coincide. This has been studied before, and we
could give some new results, but do not do so here for reasons of
space. We only mention the MMM since it comes up later in another example.
%

%
%
\subsection{Easy solutions for the process $\Vn(1) = \vtwo$}\label{sec5.1}
%
%

%
In terms of complexity, the BSDE (\ref{equ2.18}) or one of its equivalent
forms (\ref{equ3.1}), (\ref{equ4.2}), (\ref{equ4.5}) is the most difficult
one. So we focus on that equation, in the form (\ref{equ4.5}), and we try
to have a solution tuple $(Y, \varphi, W, \Ls)$ with $\varphi\equiv
0$ and
$W\equiv0$. Then (\ref{equ4.5}) simplifies to
\[
Y_t = Y_0 + \int_0^t {\lambda_s^2 \spP{Y}_s \over1+\lambda_s^2
\Delta\MM_s} \dd\MM_s + \Ls_t,
\]
which gives
$
\Delta B^Y = {\lambda^2 {}^{\p}Y \over1+\lambda^2
\Delta\MM} \Delta\MM
$.
But
$
{}^{\p}Y = Y_- + \Delta B^Y
$
by (\ref{equ2.10}), and plugging this in above and solving for $\Delta
B^Y$ allows us to get
$
{}^{\p}Y = Y_- (1+\lambda^2 \Delta\MM)
$
so that (\ref{equ4.5}) becomes
%
\begin{equation}\label{equ5.0}
Y_t = Y_0 + \int_0^t Y_{s-} \lambda_s^2 \dd\MM_s + \Ls_t,\qquad Y_T=1.
\end{equation}
This is the equation for a generalized stochastic exponential, and so
it is not surprising that we can find an explicit solution.
%
\begin{corollary}\label{Cor5.1}
Set $K:= \langle\lambda\SIdot M\rangle
$ and suppose that
\[
\E(K)_T^{-1} = c + m_T
\]
with a constant $c>0$ and a $P$-martingale $m$ which is strongly
$P$-orthogonal both to $\Mc$ and to the space of stochastic integrals
$
\{ \bar W*(\mu^S-\nu) | \bar W\in\G^2_\loc(\mu^S) \}
$.
Then the solution of (\ref{equ4.5}) is given by $\varphi\equiv0$, $W\equiv0$ and
%
\begin{eqnarray}\label{equ5.1}
Y_t
&=&
E[ \E(K)_t / \E(K)_T | \F_t]
=
\E(K)_t (c+m_t),
\nonumber\\[-8pt]\\[-8pt]
\Ls_t
&=&
\int_0^t \E(K)_{s-} \dd m_s + [ \E(K), m ]_t.
\nonumber
\end{eqnarray}
\end{corollary}
\begin{pf}
Since (\ref{equ5.0}) can be written as
$
Y = Y_0 + \int Y_- \dd K + \Ls
$,
defining $Y$ and $\Ls$ by (\ref{equ5.1}) gives by the product rule that
$(Y, \Ls)$ satisfy (\ref{equ5.0}) with $Y_T=1$, and $\Ls$ is a local
$P$-martingale like $m$ by Yoeurp's lemma. Finally, for every $\Wb
\in
\G^2_\loc(\mu^S)$, we have that
\begin{eqnarray*}
\bigl[ \Wb*(\mu^S-\nu), [\E(K),m] \bigr]
&=&
\sum\Delta\bigl( \Wb*(\mu^S-\nu) \bigr) \Delta\E(K) \Delta m\\
&=&
\Delta\E(K) \SIdot[ \Wb*(\mu^S-\nu), m]
\end{eqnarray*}
is a local $P$-martingale because $m$ is strongly $P$-orthogonal to
$\Wb
*(\mu^S-\nu)$. Hence $\Ls$ is also strongly $P$-orthogonal to $\Wb
*(\mu
^S-\nu)$, and so $(Y,0,0,\Ls)$ is a solution to (\ref{equ4.5}).
\end{pf}
%
\begin{example}\label{Exa5.2}
A special case of Corollary~\ref{Cor5.1} occurs if the
(final) \textit{mean-variance tradeoff} $\langle\lambda\SIdot
M\rangle_T$
and all the jumps $\lambda^2\Delta\MM$ are \textit{deterministic}. Then
$m\equiv0$, the solution for $Y$ is
\[
Y_t = \E( \langle\lambda\SIdot M\rangle)_t / \E( \langle\lambda
\SIdot
M\rangle)_T,
\qquad 0\le t\le T
\]
[which is adapted because $\E( \langle\lambda\SIdot M\rangle)_T$ is
deterministic], and all other quantities in the BSDEs (\ref{equ2.18}) or
(\ref{equ4.2}) or (\ref{equ4.5}) are identically 0. If $S$ or $M$ or even
only $A = \int\lambda^2\dd\MM$ is continuous, the above expression
simplifies to
\[
Y_t = e^{\langle\lambda\sSIdot M\rangle_t - \langle\lambda\sSIdot
M\rangle_T},\qquad
0\le t\le T.
\]
\end{example}

Similar results as in this section, but under more restrictive
assumptions, have been obtained by several authors. We only mention
exemplarily the work of \citet{BiaGuaPra00},
\citet{ManTev03N2} and \citet{San06}.

%
\subsection{The discrete-time case}\label{sec5.2}
%
%

%
Now we briefly look at the special case of a model in finite discrete
time $k=0,1,\ldots,T$. Our price process is given by $S =
(S_k)_{k=0,1,\ldots, T}$, and we assume as in (\ref{equ2.7}) that
%
\begin{equation}\label{equ5.2}
S = S_0 + M + \lambda\SIdot\MM
\end{equation}
with a martingale $M = (M_k)_{k=0,1,\ldots,T}$ null at 0. We assume that
$S$ is square-integrable to avoid technical complications, and we write
$\Delta_k Y:= Y_k - Y_{k-1}$ for the increments of a process $Y =
(Y_k)_{k=0,1,\ldots,T}$. The Doob decomposition
$
S = S_0 + M + A
$
is then explicitly given by
$\Delta_k A = E[\Delta_k S | \F_{k-1} ]
$,
we have
$
\Delta_k \MM= E[ (\Delta_k M)^2 | \F_{k-1} ] = \Var[\Delta_k
S
| \F_{k-1} ]
$,
and so (\ref{equ5.2}) takes the form
$
S = S_0 + M + \sum_j \lambda_j \Delta_j \MM
$
with
%
\begin{equation}\label{equ5.3}
\lambda_j = {\Delta_j A \over\Delta_j \MM} = {E[\Delta_j S |
\F
_{j-1} ] \over\Var[\Delta_j S | \F_{j-1} ]}.
\end{equation}

For the discrete-time version of the BSDE (\ref{equ2.18}), we need
$
{}^{\p}Y_j = E[ Y_j | \F_{j-1} ]
$\vspace*{2pt}
and the density $g(Y)$ of $\sdpP{[ \NY, [S] ]}$ with respect to
$\MM
$. But we have
$
[ \NY, [S] ] = \sum_j (\Delta_j \NY) (\Delta_j S)^2
$
so that
%
\begin{equation}\label{equ5.4}
g_j(Y) \Delta_j \MM= E[ (\Delta_j \NY) (\Delta_j S)^2 | \F_{j-1}].
\end{equation}
Moreover, we have
%
\begin{eqnarray}\label{equ5.5}
(1+\lambda_j^2 \Delta_j\MM) \Delta_j \MM
&=&
\Var[\Delta_j S | \F_{j-1} ] + ( E[\Delta_j S | \F_{j-1}
])^2\nonumber\\[-8pt]\\[-8pt]
&=&
E[ (\Delta_j S)^2 | \F_{j-1} ],
\nonumber
\end{eqnarray}
and the Galtchouk--Kunita--Watanabe decomposition
$
\NY= \sum_j \psi_j \Delta_j M + L
$
yields
%
\begin{eqnarray}
\label{equ5.6}
\psi_j \Delta_j \MM
&=&
\Cov( \Delta_j \NY, \Delta_j M | \F_{j-1} )
=
\Cov( \Delta_j Y, \Delta_j S | \F_{j-1} )\nonumber\\[-8pt]\\[-8pt]
&=&
\Cov( Y_j, \Delta_j S | \F_{j-1} ).\nonumber
\end{eqnarray}
Hence we get
\begin{eqnarray*}
(\psi_j + \lambda_j {}^{\p}Y_j)^2 (\Delta_j \MM
)^2
&=&
\bigl( \Cov( Y_j, \Delta_j S | \F_{j-1} ) + E[\Delta_j S |
\F
_{j-1}] E[Y_j | \F_{j-1}] \bigr)^2
\\
&=&
( E[Y_j \Delta_j S | \F_{j-1}])^2.
\end{eqnarray*}
Writing out the discrete-time analog of (\ref{equ2.18}), expanding the
ratios in the first appearing sum with $\Delta_j \MM$ and using
(\ref{equ5.3})--(\ref{equ5.6}) then yields
%
\begin{eqnarray}
\label{equ5.7}\qquad
Y_k
&=&
Y_0 + \sum_{j=1}^k {(\psi_j + \lambda_j {}^{\p}Y_j)^2 \over{}^{\p}Y_j
(1+\lambda_j^2 \Delta_j \MM) + g_j(Y)} \Delta_j \MM
+ \sum_{j=1}^k \psi_j \Delta_j M + L_k
\nonumber\\
&=&
Y_0 + \sum_{j=1}^k { (E[Y_j \Delta_j S | \F_{j-1}])^2 \over
E[Y_j
| \F_{j-1}] E[(\Delta_j S)^2 | \F_{j-1}] + E[ (\Delta_j \NY)
(\Delta_j S)^2 | \F_{j-1} ] }
\\
&&{}
+ \sum_{j=1}^k \psi_j \Delta_j M + L_k,\qquad
Y_T=1.
\nonumber
\end{eqnarray}
But
$
Y_j = Y_{j-1} + \Delta_j \NY+ \Delta_j B^Y
$
gives
\[
E[Y_j | \F_{j-1}] = Y_{j-1} + \Delta_j B^Y = \NY_{j-1} + B^Y_j,
\]
and the denominator in the third sum in (\ref{equ5.7}) therefore equals
\[
E[ (\NY_{j-1} + B^Y_j + \Delta_j \NY) (\Delta_j S)^2 | \F
_{j-1} ]
=
E[ Y_j (\Delta_j S)^2 | \F_{j-1} ].
\]
Passing to increments and taking conditional expectations to make the
martingale increments vanish, equation (\ref{equ5.7}) thus can be
written as
\[
Y_{k-1}
=
E[ Y_k - \Delta_k Y | \F_{k-1} ]
=
E[ Y_k | \F_{k-1} ]
- { (E[Y_k \Delta_k S | \F_{k-1} ])^2 \over E[Y_k (\Delta_k S)^2
| \F_{k-1}] },\qquad
Y_T=1.
\]
This is exactly the recursive relation derived in equation (3.1) in
Theorem 1 of \citet{Gug03}; see also equation (3.36) in
\citet{CerKal07}. Under more
restrictive assumptions, analogous equations have also been obtained in
equation (5) in Theorem 2 of \citet{autokey5} or in equation
(2.19) in Theorem 1
of \citet{BerKogLo01}.

\subsection{\texorpdfstring{On the relation to Arai (\citeyear{Ara05})}{On the relation to Arai (2005)}}\label{sec5.3}

Our final example serves to illustrate the relations between our work
and that of \citet{Ara05}, whose assumptions are rather similar to ours. More
precisely, \citet{Ara05} assumes that $S$ (which he calls $X$) is locally
bounded, and that the VOMM\vspace*{-1pt} $\Qt$ exists in $\PPes$ and satisfies the
reverse H\"older inequality $R_2(P)$ and a condition on the jumps of
$Z^{\Qt}$. This implies of course $S\in\SS^2_\loc(P)$ and $\PPes
\ne
\varnothing$. \citet{Ara05} does not use BSDEs, but works with a change of
numeraire as in \citet{GouLauPha98}. His numeraire is
$E_{\Qt}[ Z^{\Qt}_T | \F
_{\fatdot}]$, and to ensure that this is positive, the existence of the
VOMM $\Qt$ in $\PPes$ is needed. The example below illustrates that our
assumptions are strictly weaker than those of \citet{Ara05}.
%
\begin{example}\label{Exa5.3}
We start with two independent simple Poisson processes
$N^{(\pm)}$ with the same intensity $\alpha>0$ and define $n^\pm_t:=
N^{(\pm)}_t - \alpha t$, $0\le t\le T$. We then set
\[
dS_t = S_{t-} (\gamma_+ \dd n^+_t - \gamma_- \dd n^-_t + \delta\dd t)
=: S_{t-} \dd R_t,
\]
so that $S$ is clearly locally bounded, hence in $\SS^2_\loc(P)$, and
even quasi-left-continuous. We claim that we can choose the parameters
$\alpha, \gamma_+, \gamma_-, \delta$ such that:

\begin{longlist}[(2)]
\item[(1)]
$\PPes\ne\varnothing$;

\item[(2)] the variance-optimal signed martingale measure $\Qt\in
\PPss$
coincides with the minimal signed martingale measure $\Qh$, but is not
in $\PPes$, which means in our terminology and that of \citet{Ara05}
that the
VOMM does not exist.

Let us first argue (2). Because $dM_t = S_{t-} (\gamma_+ \dd n^+_t -
\gamma
_-\dd n^-_t)$ implies that
$d\MM_t = S_{t-}^2 (\gamma_+^2+\gamma_-^2)\alpha\dd t$ and we have
$dA_t = S_{t-} \delta\dd t$, we obtain
\[
\lambda\SIdot M = {\delta\over\alpha(\gamma_+^2 + \gamma_-^2)}
(\gamma
_+ n^+ - \gamma_- n^-).
\]
So as soon as we have
%
\begin{equation}\label{equ5.8}
{\delta\gamma_+ \over\alpha(\gamma_+^2 + \gamma_-^2)} > 1,
\end{equation}
we get $- \lambda\Delta M < -1$ at jumps of $N^{(+)}$ so that $\Zh=
\E
(-\lambda\SIdot M)$ also takes negative values. Because the
mean-variance tradeoff process\vadjust{\goodbreak}
$
\langle\lambda\SIdot M\rangle_t = \break{\delta^2 \over\alpha(\gamma
_+^2 +
\gamma_-^2)} t
$,
$0\le t\le T$, is deterministic, the signed MMM $\Qh$ is
variance-optimal by\vspace*{1pt} Theorem~8 of \citet{Sch95}. Moreover, $\Zh$ is
clearly in
$\M^2(P)$ and so $\Qt=\Qh$ is in $\PPss$, but not in $\PPes$. This
gives (2).

To construct an element of $\PPes$, start with
$
Z:= \E(L):= \E(\beta_1 n^+ + \beta_2 n^-)
$,
which is clearly in $\M^2(P)$. To ensure that $Z>0$, we need $\beta
_1>-1$ and $\beta_2>-1$. Next, the product $ZS$ is by It\^o's
formula seen to be a local $P$-martingale if and only if
$
\delta\dd t + d\langle L, R\rangle_t \equiv0
$,
which translates into the condition
$
\delta= (\beta_2 \gamma_- - \beta_1 \gamma_+) \alpha
$.
This allows us to rewrite (\ref{equ5.8}) as
\[
{\gamma_+^2 + \gamma_-^2 \over\gamma_+} < {\delta\over\alpha} =
\beta
_2 \gamma_- - \beta_1 \gamma_+,
\]
and if we choose $\gamma_+ = \gamma_- = \gamma$, this boils down to
$
\beta_2 - \beta_1 > 2
$
and
$
{\delta\over\alpha} = (\beta_2 - \beta_1) \gamma
$.
By the Bayes rule, $S$ is then a local $Q$-martingale under $Q\approx
P$ with $dQ = Z_T \dd P$.

If we now choose $\varepsilon>0$ and
$
\beta_1=\beta> -1
$,
$
\beta_2 = \beta+2+\varepsilon
$,
$
\alpha=1
$,
$
\delta=(2+\varepsilon)\gamma
$,
one readily verifies that all conditions above are satisfied; hence
$\PPes\ne\varnothing$ since it contains $Q$. If we take $\gamma\in(0,1)$,
we even keep $S>0$ since $\Delta R > -1$.
\end{longlist}
\end{example}
\begin{Rem*}
By its construction, the minimal martingale density $\Zh$
is always based on $-\lambda\SIdot M$. With our above choice of model
parameters $\gamma_+ = \gamma_- = \gamma$, this is symmetric in $n^+$
and $-n^-$ and therefore risks getting negative jumps rather easily. In
contrast, writing
\[
L = \beta n^+ + (\beta+2+\varepsilon) n^- = - \lambda\SIdot M + \Lt
\]
with $\Lt= (\beta+1+{\varepsilon\over2}) n^+ + (\beta
+1+{\varepsilon\over2})
n^-$ shows that it can be very beneficial to have some extra freedom
when choosing an ELMM or a martingale density. This is quite analogous
to the well-known counterexample in \citet{DelSch98}.
\end{Rem*}



\printaddresses

\end{document}